\documentclass{article}  
\usepackage[english]{babel}
\usepackage{amsmath,amsthm,amssymb,amscd}
\usepackage[cp1251]{inputenc}
\input{amssym.def}
\newtheorem{thm}{Theorem}[section]
\newtheorem{prop}[thm]{Proposition}

\newtheorem{defn}[thm]{Definition}
\def\proof{\par\noindent {\sl Proof. }}
\newcommand{\rf}{\,\mathsf r^f\,} 
\newcommand{\ro}{\,\mathsf r^{f^\circ}\,}
\newcommand{\rk}{\,\mathsf r\,}
\def\fin{\hfill$\Box$\medskip}
\title{On certain constructive predicate calculus\thanks{The work was supported by RFBR grant 20-01-00670.}}
\author{V.\,E.~Plisko\thanks{Faculty of Mechanics and Mathematics of Lomonosov Moscow State University; veplisko@yandex.ru}}
\date{}
\begin{document}
\maketitle
\begin{abstract}
Constructive arithmetic, or the Markov arithmetic MA, is obtained from intuitionistic arithmetic HA by adding the following two principles: the 
Markov principle M which distinguishes constructivism from intuitionism, and the so-called extended Church thesis ECT, which distinguishes 
constructive semantics from classical semantics.
The first pinciple is expressed by a predicate formula, but ECT is formulated as a scheme in the arithmetical language.
Some earlier results of the authors make possible to replace the scheme ECT by a pure predicate formula.
This gives a predicate calculus MQC which can serve as a logical basis for constructive arithmetic.
Namely, arithmetical theory based on MQC and Peano axioms proves all formulas deducible in MA.
Note that MQC is not an intermediate calculus between intuitionistic and classical logics: it proves some formulas that are not deducible in the
classical predicate calculus.

Bibliography: 12 items.
\end{abstract}
\section{Introduction} \label{sv}

Intuitionism connects the truth of a statement with its provability, thus a statement is considered as true if there is a proof of it.
In the context of such an interpretation of truth, logical operations receive a peculiar in\-ter\-pre\-ta\-tion.
Starting with the work of Kolmogorov~\cite{Kolm3}, much attention is paid to logical systems that are correct from the point of view of 
intuitionism.
The intuitionistic predicate calculus $\mathsf{IQC}$ was developed, and intuitionistic arithmetic $\mathsf{HA}$ is based on $\mathsf{IQC}$.
Then Kleene~\cite{Kleene1} proposed  recursive realizability as an interpretation of specific intuitionistic concepts based on the theory of 
recursive functions.
Re\-cur\-si\-ve realizability can be considered as a kind of semantics of mathematical statements.
This semantics underlies the constructive approach in mathematics, the systematic development of which was carried out by Markov and his students 
and followers.
In the course of the study of constructive semantics, some logical and mathematical principles characteristic for constructivism were identified 
and formulated: the Markov principle, or the principle of constructive selection (see~\cite{Mark2}), which distinguishes constructivism from 
intuitionism, and the so-called extended Church thesis, which distinguishes constructive semantics from classical one.
These two principles, the Markov principle and the extended Church thesis, when added to the axiomatics of intuitionistic arithmetic 
$\mathsf{HA}$ give the theory $\mathsf{MA}$, sometimes called the Markov arithmetic or traditional constructivism.

Later, the author proved that the predicate logic of recursive realizability is not arithmetical (see~\cite{Plisko13}).
Moreover, this semantics of predicate formulas is largely occasional: one can expand the arithmetical language and define the concept of 
realizability for it in such a way that the class of realizable formulas will narrow.
This procedure can be continued by transfinite induction to any constructive ordinal (see~\cite{Plisko15}).
As a result, the author proposed the concept of absolute realizability of predicate formulas, independent of the language in which predicates 
interpreting predicate variables are formulated (see~\cite{Plisko16}).

A by-product of the author's research was the so-called scheme theorem which makes it possible to translate predicate schemes over the 
language of arithmetic into pure predicate formulas while preserving constructive validity.
In particular, the predicate scheme expressing a variant of the extended Church thesis can also be replaced by a predicate formula.
The result is a predicate calculus $\mathsf{MQC}$ as an extension of $\mathsf{IQC}$ by the Markov principle and a pure predicate 
translation of the extended Church thesis.
There is every reason to call the proposed calculus a constructive predicate calculus.

Section~\ref{ska} describes the constructive semantics of the formal arithmetical language and its pure predicate variant, discusses some arithmetical theories sound relative to this semantics, and also establishes a number of facts used later.
Section~\ref{sklp} describes the constructive semantics of predicate formulas based on the concept of absolute realizability.
Section~\ref{stos} outlines the main technical result, namely the scheme theorem.
Finally, Section~\ref{skip} presents the calculus $\mathsf{MQC}$, in which only absolutely realizable predicate formulas are derived.
It is proved that this calculus can serve as a logical basis for constructive Markov arithmetic.

\section{Constructive arithmetic} \label{ska}

\subsection{Intuitionistic arithmetical theories}

The signature of {\it the  language of formal arithmetic} $LA$ consists of the individual constant $0$, the functional symbol $s$ for the unary operation $x\mapsto x+1$, the functional symbols $+$ and $\cdot$ for addition and multiplication, and the 
predicate equality symbol ~$=$.
The alphabet of the language also contains propositional connectives $\&$, $\lor$, $\to$, $\neg$, quantifier symbols $\forall$ and~$\exists$.
We will not distinguish between {\it natural numbers} $0,1,2,\dots$ and the {\it terms} $0,s0,ss0,\dots$ representing them.
By $\bot$ we denote the $LA$-formula $0=1$.
The expression $\Phi\leftrightarrow\Psi$, where $\Phi$ and $\Psi$ are formulas, will be considered as an abbreviated notation for $(\Phi\to\Psi)\,\&\,(\Psi\to\Phi)$.
Sometimes instead of $\forall x_1\ldots\forall x_n\,\Phi$, we will write $\forall x_1,\ldots,x_n\,\Phi$ or even $\forall\vec x\,\Phi$, where
$\vec x$ is the list of variables $x_1,\ldots,x_n$.
If the notation $\Phi(x_1,\ldots,x_n)$ is used for a formula $\Phi$, it means that $\Phi$ does not contain free variables other than
$x_1,\ldots,x_n$.
Then $\Phi(t_1,\ldots,t_n)$ denotes the result of substituting in $\Phi$ terms $t_1,\ldots,t_n$ for free occurrences of variables
$x_1,\ldots,x_n$ respectively.
In this case, the bound variables in $\Phi$ are renamed  in such a way that the substitution becomes free in the sense 
of~\cite[\S~18]{Kleene3}.
If $t_1$ and $t_2$ are arbitrary terms, then $t_1\le t_2$ will denote $\exists v\,(t_1+v=t_2)$, where $v$ is a variable not in $t_1$ and $t_2$.
{\it Bounded quantifiers} $(\forall x\le t)$ and $(\exists x\le t)$, where $t$ is a term that does not contain the variable $x$, are sometimes 
used to abbreviate arithmetical formulas.
Namely, $(\forall x\le t)\,\Phi$ and $(\exists x\le t)\,\Phi$ are considered as abbreviated forms of $\forall x\,(x\le t\to\Phi)$ and
$\exists x\,(x\le t\,\&\,\Phi )$ respectively.

Arithmetical {\it statements} are closed $LA$-formulas.
The classical {\it truth} of an arithmetical statement is its truth in the {\it standard interpretation} $\mathfrak N$ of the language $LA$.

{\it Intuitionistic arithmetic} $\mathsf{HA}$ is a theory based on {\it intuitionistic predicate calculus} $\mathsf{IQC}$ and Peano axioms in the 
language $LA$ including the axiom scheme of induction.
{\it Intuitionistic Robinson arithmetic} $\mathsf{IRA}$ is a theory defined by a finite number of non-logical axioms in the language $LA$ and
obtained by replacing in $\mathsf{HA}$ the axiom scheme of induction with the axiom $\forall x\,(x=0\lor\exists y\,(x=sy))$.

In what follows, when proving the deducibility of formulas in the theory $\mathsf{HA}$ or its extensions, we will use the following technique
of {\it proof by induction}: by proving $\Gamma\vdash\Phi(0)$ and $\Gamma\vdash\Phi(x)\to\Psi(sx)$, we conclude that $\Gamma\vdash\Phi(x)$.

By $LA^f$ we denote the extension of the language $LA$ by introducing symbols for all primitive recursive functions, and we assume that each
such symbol encodes a way to obtain this primitive recursive function from the basic functions
$$o(x)=0,\,s(x)=x+1,\,I^n_k(x_1,\ldots,x_n)=x_k (n\ge 1, 1\le k\le n)$$
using substitution and recursion.
The theory $\mathsf{HA}^{\mathsf f}$ is an extension of the theory $\mathsf{HA}$ in the language $LA^f$ by adding defining equalities for 
all additional functional symbols.

We will need the ability to encode lists of natural numbers with natural numbers.
To do this, we will introduce suitable primitive recursive functions into consideration.
A binary primitive recursive function $\mathsf c$, defined by the formula $\mathsf c(x,y)=\frac{(x+y)(x+y+1)}{2}+x$, sets the one-to-one numbering 
of all pairs of natural numbers, and there are primitive recursive functions $\mathsf l$ and $\mathsf r$ such that 
$\mathsf c(\mathsf l(x),\mathsf r(x))=x$, $\mathsf l(\mathsf c(x,y))=x$, $\mathsf r(\mathsf c(x,y))=y$ are derivable in $\mathsf{HA}$.
The binary function $\mathsf l_*$ is defined by the following recursive scheme (and therefore is primitive recursive): $\mathsf l_*(x,0)=x$,
$\mathsf l_*(x,sy)=\mathsf l(\mathsf l_*(x,y))$.
The binary primitive recursive function $\mathsf g$ is defined as follows: $\mathsf g(x,y)=\mathsf r(\mathsf l_*(x,y))$.
It is not difficult to make sure that the function $\mathsf g$ performs one-to-one numbering of all total functions of the type 
$\mathbb N\to\mathbb N$ that take a value other than 0 only in a finite number of points.
The ternary function $\mathsf h$ is defined as follows: $\mathsf h(x,y,0)=\mathsf c(\mathsf l(x),y)$,
$\mathsf h(x,y,sz)=\mathsf c(\mathsf h(\mathsf l(x),y,z),\mathsf r(x))$.
It can be shown that $\mathsf h$ is a primitive recursive function and the following formulas are derived in $\mathsf{HA}^f$:

\begin{equation} \label{f7} 
\neg(v=z)\to\mathsf g(\mathsf h(x,y,z),v)=\mathsf g(x,v),
\end{equation}

\begin{equation} \label{f83} 
\mathsf g(\mathsf h(x,y,z),z)=y.
\end{equation}

Denote by $B(u,v)$ the formula $\mathsf g(u,v)=1$.
         
\begin{prop} \label{p13}
For any $LA^f$-formula $\Phi(x)$, the formula 
$$\forall x\,\neg\neg\exists y\forall z\,(z\le x\to(B(y,z)\leftrightarrow\Phi(z)))$$
is deducible in $\mathsf{HA}^f$. 
\end{prop}

\proof
Using the proof by induction, we establish that the formula 
\begin{equation} \label{f15}
\neg\neg\exists y\forall z\,(z\le x\to(B(y,z)\leftrightarrow\Phi(z)))
\end{equation}
is deducible in $\mathsf{HA}^f$.
First we show that the formula 
$$\neg\neg\exists y\forall z\,(z\le 0\to(B(y,z)\leftrightarrow\Phi(z)))$$
is deducible in $\mathsf{HA}^f$.
Obviously, it is enough to make sure that the formula 
\begin{equation} \label{f8}
\exists y\forall z\,(z\le 0\to(B(y,z)\leftrightarrow\Phi(z))) 
\end{equation}
is deducible from the hypothesis $\Phi(0)\lor\neg\Phi(0)$.
Let's use the rule of case analysis, or elimination of disjunction, namely, we show that~(\ref{f8}) is derived from each of the 
hypotheses $\Phi(0)$ and $\neg\Phi(0)$.

Fix the hypothesis $\Phi(0)$, put $a=\mathsf c(0,1)$ and prove the deducibility of $B(a,z)\leftrightarrow\Phi(0)$ from the hypothesis 
$z\le 0$, which is equivalent to $z=0$.
It is enough to prove the deducibility in $\mathsf{HA}^f$ of $\mathsf g(a,0)=1$, and this is justified
by the following calculations:
$$\mathsf g(\mathsf c(0,s0),0)=\mathsf r(\mathsf l_*(\mathsf c(0,s0),0))=\mathsf r(\mathsf c(0,s0))=s0.$$

Now fix the hypothesis $\neg\Phi(0)$, put $a=\mathsf c(0,0)$ and prove the deducibility of $B(a,z)\leftrightarrow\Phi(0)$ from the 
hypothesis $z=0$.
Obviously, it is enough to prove the deducibility in $\mathsf{HA}^f$ of $\neg(\mathsf g(a,0)=1)$, and this follows from the 
deducibility of $\mathsf g(a,0)=0$, which is justified by the following calculations:
$$\mathsf g(\mathsf c(0,0),0)=\mathsf r(\mathsf l_*(\mathsf c(0,0),0))=\mathsf r(\mathsf c(0,0))=0.$$

Now we show that $\neg\neg\exists y\forall z\,(z\le sx\to(B(y,z)\leftrightarrow\Phi(x)))$ is derivable in $\mathsf{HA}^f$ 
from the hypothesis~(\ref{f15}).
Obviously, it is enough to make sure that the formula 
\begin{equation} \label{f11}
\exists y\forall z\,(z\le sx\to(B(y,z)\leftrightarrow\Phi(z))) 
\end{equation}
is deducible from the hypotheses
\begin{equation} \label{f115}
\forall z\,(z\le x\to(B(y,z)\leftrightarrow\Phi(x)))
\end{equation}
and $\Phi(sx)\lor\neg\Phi(sx)$.
We show that~(\ref{f11}) is derived from~(\ref{f115}) and each of the hypotheses $\Phi(sx)$ and $\neg\Phi(sx)$.

In the case of the hypothesis $\Phi(sx)$ we denote by $t$ the term $\mathsf h(y,sx,s0)$ and prove that
$B(t,z)\leftrightarrow\Phi(z)$ is derived from the hypothesis $z\le sx$.
The formula $(z\le x)\lor(z=sx)$ is derived from this hypothesis, therefore, each of the hypotheses $z\le x$ and $z=sx$ can be considered 
separately.
The formula $B(y,z)\leftrightarrow\Phi(z)$ is derived from the hypotheses $z\le x$ and~(\ref{f115}).
It remains to be noted that the deducibility of~(\ref{f7}) implies the deducibility of the equality $g(t,z)=g(y,z)$.
In the case of the hypothesis $z=sx$, it is sufficient to prove the deducibility of the equality $g(t,sx)=s0$, and this immediately follows from 
the deducibility of~(\ref{f83}).

The hypothesis $\neg\Phi(x)$ is treated in exactly the same way if we take the term $\mathsf h(y,sx,0)$ as~$t$.~\fin

$\Delta_0$-{\it formulas} of the language $LA$ are defined inductively as follows:
1)~every atomic formula is a $\Delta_0$-formula;
2)~if $\Phi$ and $\Psi$ are $\Delta_0$-formulas, then $\neg\Phi$, $(\Phi\vee\Psi)$, $(\Phi\,\&\,\Psi)$, $(\Phi\to\Psi)$ are
$\Delta_0$-formulas;
3)~if $x$ is a variable, $t$ is a term that does not contain $x$, $\Phi$ is a $\Delta_0$-formula, then $(\exists x\le t)\,\Phi$ and
$(\forall x\le t)\,\Phi$ are $\Delta_0$-formulas.

Arithmetical $\Sigma$-{\it formulas} are defined inductively as follows:
1)~every $\Delta_0$-formula is a $\Sigma$-formula;
2)~if $\Phi$ and $\Psi$ are $\Sigma$-formulas, then $(\Phi\vee\Psi)$ and$(\Phi\,\&\,\Psi)$ are also $\Sigma$-formulas;
3)~if $x$ is a variable, $t$ is a term that does not contain occurrences of $x$, and $\Phi$ is a $\Sigma$-formula, then $(\exists x\le t)\,\Phi$, 
$(\forall x\le t)\,\Phi$, $\exists x\,\Phi$ are $\Sigma$-formulas.

As it follows from~\cite[\S\S~41,49,74]{Kleene3}, the theory $\mathsf{HA}^{\mathsf f}$ is a conservative extension of the theory $\mathsf{HA}$ since all additional function symbols can be eliminated: for every primitive recursive function $f:\mathbb N^n\to\mathbb N$ a $\Sigma$-formula 
$\Phi(x_1,\dots,x_n,y)$ is constructed in such a way that for any natural $k_1,\ldots,k_n,\ell$, if $f(k_1,\ldots,k_n)=\ell$, then $\Phi(k_1,\ldots,k_n,\ell)$ is derived in $\mathsf{IRA}$, and otherwise $\neg\Phi(k_1,\ldots,k_n,\ell)$ is derived, and at the same 
time the formula 
\begin{equation} \label{f10}
\forall x_1,\ldots,x_n\exists!y\,\Phi(x_1,\dots,x_n,y)
\end{equation}
is derived in $\mathsf{HA}$.
In the theory $\mathsf{HA}^{\mathsf f}$, the defining equalities for the function $f$ can be replaced by the axioms~(\ref{f10}) and
$\forall x_1,\ldots,x_n\,\Phi(x_1,\dots,x_n,f(x_1,\dots,x_n))$.
Iin the theory $\mathsf{HA}^{\mathsf f}$ the formula $f(x_1,\dots,x_n)=y$ is provably equivalent to 
$\Phi(x_1,\dots,x_n,y)$.
This allows us to define a translation from the language $LA^f$ to the language $LA$, which ensures the conservativeness of the theory 
$\mathsf{HA}^{\mathsf f}$ over $\mathsf{HA}$.
Namely, for each $LA^f$-formula $\Phi$, the $LA$-formula $\Phi'$ is constructed as in~\cite[\S~74, lemma~29]{Kleene3}, so the following
statement holds:

\begin{prop}  \label{p19} 
$\Gamma\vdash_{\mathsf{HA}^{\mathsf f}}\Phi$ if and only if $\Gamma'\vdash_{\mathsf{HA}}\Phi'$, where $\Gamma'$ is the list of $LA$-formulas
$\Psi_1',\ldots,\Psi_n'$ if $\Gamma$ is the list of $LA^f$-formulas $\Psi_1,\ldots,\Psi_n$.
\end{prop}

Proposition~\ref{p19} allows us to take some liberties in expressions: speaking further about the deducibility of a $LA^f$-formula $\Phi$ in the 
theory $\mathsf{HA}$, we will keep in mind the deducibility of $LA$-formula $\Phi'$.
In particular, we note the following fact.

\begin{prop}  \label{p20}
For any atomic formula $\Phi(\vec x)$ of the language $LA^f$, where $\vec x$ is a list of variables $x_1,\ldots,x_n$, the formula 
$\forall\vec x\,(\Phi(\vec x)\lor\neg\Phi(\vec x))$ is derived in $\mathsf{HA}$.
\end{prop}

\proof
This follows from the fact that $\mathsf{HA}$ is a theory with decidable equality: it derives $\forall x,y\,(x=y\lor\neg(x=y))$.~\fin

Proposition~\ref{p20} implies
\begin{equation} \label{f21}
\vdash_{\mathsf{HA}}\forall x,y\,(B(x,y)\lor\neg B(x,y)),
\end{equation}
where $B(x,y)$ is the formula expressing the predicate $\mathsf g(x,y)=1$.

***

\subsection{A pure predicate arithmetical language} \label{sscpvya}

Let's apply to the language $LA$ the procedure for eliminating functional symbols described, for example, in~\cite[\S~74]{Kleene3}. 
The predicate symbol of equality $=$ is also replaced by the binary predicate symbol~$E$.
As a result, we get an elementary language $Ar$, the signature of which consists of an unary predicate symbol $Z$, binary predicate symbols $E$ 
and $S$, and ternary predicate symbols $A$ and~$M$.
In this case, each $LA$-formula $\Phi$ is mapped to its {\it predicate form} $\Phi'$, as in~\cite[Lemma~29]{Kleene3}.
(Note that each $Ar$-formula can be considered as a predicate form of some $LA$-formula.)

In the standard interpretation, the symbols $Z,E,S,A,M$ have the following meaning: $Z(x)$ means $x=0$, $E(x,y)$ means $x=y$, 
$S(x,y)$ means $x+1=y$, $A(x,y,z)$ means $x+y=z$, $M(x,y,z)$ means $x\cdot y=z$.

Along with the language $Ar$, we will consider the elementary language $Ar^*$, the signature of which is obtained by adding the individual
constants $0,1,2,\ldots$ for all natural numbers to the signature of the language $Ar$.
Also, we denote the $Ar^*$-formula $E(0,1)$ by $\bot$.

We will say that a closed $Ar^*$-formula $\Phi$ {\it is true} if $\mathfrak N\models\Phi$ in the usual classical sense.
Closed formulas of the language $Ar^*$ will also be called arithmetical statements.

There is a natural translation from the language $Ar^*$ to the language~$LA$.
It consists of constructing for each $Ar^*$-formula $\Phi$ its {\it functional notation}, namely, $LA$-formula $\Phi^\circ$ obtained by replacing
each constant $n$ in it with $LA$-term $n$ and atomic formulas $Z(t)$, $E(t,u)$, $S(t,u)$, $A(t,u,v)$ and $M(t,u,v)$, where $t,u,v$ are constants 
or variables, for $t=0$, $t=u$, $st=u$, $t+u=v$, and $t\cdot u=v$, respectively.

The elimination of functional symbols in a theory $\mathsf T$ in the language $LA$ consists of constructing its predicate 
version, namely, a theory $\mathsf T^*$ in the language $Ar$ such that an $LA$-formula $\Phi$ is deducible in $\mathsf T$ if and only 
if its predicate form $\Phi^*$ is deducible in $\mathsf T^*$.
The procedure for constructing a predicate version of the theory $\mathsf T$ is described in~\cite[\S~74]{Kleene3}. 
In particular, applying this procedure to the theory $\mathsf{IRA}$ gives the theory $\mathsf{IRA}^*$ with the following axioms in the 
language~$Ar$:

\noindent
$A_1$. $\forall x,y,z\,(S(x,z)\,\&\,S(y,z)\to E(x,y))$\\ 
$A_2$. $\forall x,y\,(S(x,y)\to\neg Z(y))$\\
$A_3$. $\forall x,y,z\,(E(x,y)\,\&\,E(x,z)\to E(y,z))$\\ 
$A_4$. $\forall x,y,z\,(E(x,y)\,\&\,S(x,z)\to S(y,z))$\\ 
$A_5$. $\forall x,y\,(Z(y)\to A(x,y,x))$\\ 
$A_6$. $\forall x,y,z,u,v\,(S(y,z)\,\&\,A(x,y,u)\,\&\,S(u,v)\to A(x,z,v))$\\ 
$A_7$. $\forall x,y\,(Z(y)\to M(x,y,y))$\\
$A_8$. $\forall x,y,z,u,v\,(S(y,z)\,\&\,M(x,y,u)\,\&\,A(u,x,v)\to M(x,z,v))$\\
$A_9$. $\forall x_1,x_2,y,z\,(E(x_1,x_2)\,\&\,A(x_1,y,z)\to A(x_2,y,z))$\\
$A_{10}$. $\forall x,y_1,y_2,z\,(E(y_1,y_2)\,\&\,A(x,y_1,z)\to A_2,y_2,z))$\\
$A_{11}$. $\forall x_1,x_2,y,z\,(E(x_1,x_2)\,\&\,M(x_1,y,z)\to M(x_2,y,z))$\\
$A_{12}$. $\forall x,y_1,y_2,z\,(E(y_1,y_2)\,\&\,M(x,y_1,z)\to M_2,y_2,z))$\\ 
$A_{13}$. $\forall x\,(Z(x)\lor\exists y\, S(y,x)))$\\
$A_{14}$. $\exists x\,Z(x)$\\
$A_{15}$. $\forall x,y\,(Z(x)\,\&\,Z(y)\to E(x,y))$\\
$A_{16}$. $\forall x,y\,(E(x,y)\,\&\,Z(x)\to Z(y))$\\
$A_{17}$. $\forall x\,\exists y\,S(x,y)$\\ 
$A_{18}$. $\forall x,y,z\,(S(x,y)\,\&\,S(x,z))\to E(y,z))$\\
$A_{19}$. $\forall x,y,z\,(S(x,y)\,\&\,E(y,z))\to S(x,z))$\\ 
$A_{20}$. $\forall x,y\,\exists z\,A(x,y,z)$\\
$A_{21}$. $\forall x,y,z_1,z_2\,(A(x,y,z_1)\,\&\,A(x,y,z_2)\to E(z_1,z_2))$\\
$A_{22}$. $\forall x,y,z_1,z_2\,(A(x,y,z_1)\,\&\,E(z_1,z_2)\to A(x,y,z_2))$\\ 
$A_{23}$. $\forall x,y\,\exists z\,M(x,y,z)$\\
$A_{24}$. $\forall x,y,z_1,z_2\,(M(x,y,z_1)\,\&\,M(x,y,z_2)\to E(z_1,z_2))$\\
$A_{25}$. $\forall x,y,z_1,z_2\,(M(x,y,z_1)\,\&\,E(z_1,z_2)\to M(x,y,z_2))$

The connection between the theories $\mathsf{IRA}$ and $\mathsf{IRA}^*$ consists, in particular, of the fact that an $LA$-formula $\Phi$ is 
deducible in $\mathsf{IRA}$ if and only if its predicate form $\Phi^*$ is deducible in $\mathsf{IRA}^*$.

Denote by $x\le y$ the $Ar$-formula $\exists z\,A(x,z,y)$.
Now the notions of $\Delta_0$- and $\Sigma$-formulas of the language $Ar$ are defined literally in the same way as the corresponding notions
for the language $LA$.
In this case, the predicate form of a $\Sigma$-formula of the language $LA$ turns out to be a $\Sigma$-formula of the language $Ar$, and the 
functional notation of a $\Sigma$-formula of the language $Ar^*$ is a $\Sigma$-formula of the language $LA$, thus the notion of a $\Sigma$-formula is invariant in this sense.

The G\"odel numbering of the partial recursive functions is described in~\cite[\S~65]{Kleene3}. 
A unary partial recursive function with the G\"odel number $x$ will be denoted by $\{x\}$.
Every natural number is a G\"odel number of a unary partial recursive function.
$!\{x\}(y)$ will mean that the value $\{x\}(y)$ is defined.

The recursively enumerable predicate $\{x\}(y)=z$ is defined in $\mathfrak N$ by a $\Sigma$-formula $G(x,y,z)$ of the language $Ar$ such that
\begin{equation} \label{f33}
\vdash_{\mathsf{HA}}\forall x,y,z_1,z_2\,(G(x,y,z_1)\,\&\,G(x,y,z_2)\to z_1=z_2).
\end{equation}

Let $H(x,y)$ be the $\Sigma$-formula $\exists z\,(Z(z)\,\&\,G(x,y,z)).$
This formula defines in $\mathfrak N$ the predicate $\{x\}(y)=0$.

For the predicate form of the formula $B(v,x)$ defined above, we keep the same notation.

Consider the following $Ar$ formulas:

\noindent
$\begin{array}{ll}
A_{26}. & \forall x,y,z_1,z_2\,(G(x,y,z_1)\,\&\,G(x,y,z_2)\to E(z_1,z_2))\\
A_{27}. & \forall y,z\,\neg\neg\exists v\forall x\,(x\le z\to(B(v,x)\leftrightarrow H(y,x)))\\
A_{28}. & \forall x\forall y\,(B(x,y))\lor\neg B(x,y))
\end{array}$

Denote by $Q$ the conjunction of the formulas $A_1$-$A_{28}$.

\subsection{Recursive realizability}

Let $\pi_0,\pi_1,\pi_2,\ldots$ be consecutive primes, i.e. $\pi_0=2$, $\pi_1=3$, $\pi_2=5$, etc.
$(a)_i$ will denote the exponent of $\pi_i$ in the decomposition of the number $a$ into prime factors.

The following notion called {\it recursive realizability} is introduced by Kleene~\cite{Kleene1}.

\begin{defn} \label{d3} \rm
The relation $e\rk\Phi$, where $e\in\mathbb N$, $\Phi$ is a statement of the language $LA$ or $Ar^*$, is defined by induction on the number
of logical symbols in~$\Phi$.
\begin{itemize}
\item
$e\rk\Phi\rightleftharpoons e=0$ and $\Phi$ is true if $\Phi$ is an atomic statement.
\item
$e\rk(\Phi\,\&\,\Psi)\rightleftharpoons(e)_0\rk\Phi$ and $(e)_1\rk\Psi$.
\item
$e\rk(\Phi\lor\Psi)\rightleftharpoons(e)_0=0$ and $(e)_1\rk\Phi$ or $(e)_0=1$ and $(e)_1\rk\Psi$.
\item
$e\rk(\Phi\to\Psi)\rightleftharpoons\forall a\,(a\rk\Phi\Rightarrow!\{e\}(a)$ and $\{e\}(a)\rk\Psi)$.
\item
$e\rk\neg\Phi\rightleftharpoons e\rk(\Phi\to\bot)$.
\item
$e\rk\exists x\,\Phi(x)\rightleftharpoons(e)_1\rk\Phi((e)_0)$.
\item
$e\rk\forall x\,\Phi(x)\rightleftharpoons\forall a\,(!\{e\}(a)$ and $\{e\}(a)\rk\Phi(a))$.
\end{itemize}
\end{defn}

If $e\rk\Phi$ holds, then we say that the natural number $e$ {\it realizes} the statement $\Phi$, or is {\it a realization}
of the statement~$\Phi$.
A statement $\Phi$ is called {\it realizable} if there is a natural number $e$ that realizes~$\Phi$.

The following statement immediately follows from Definition~\ref{d3}.

\begin{prop} \label{p5} 
A statement of the form $\neg\Phi$ is realizable if and only if $0\rk\neg\Phi$.
\end{prop}

\begin{prop} \label{p69}
For any arithmetical $\Sigma$-formula $\Psi(\vec x)$, where $\vec x$ is a list of variables $x_1,\ldots,x_k$,
there is a $k$-place partial recursive function $\alpha_\Psi$ such that for any list of natural numbers $\vec m=m_1,\ldots,m_k$,
if the statement $\Psi(\vec m)$ is true, then $!\alpha_\Psi(\vec m)$ and $\alpha_\Psi(\vec m)\rk\Psi(\vec m)$.
\end{prop}

\proof
The function $\alpha_\Psi$ is defined by induction on the number of logical symbols in~$\Psi$ in accordance with Definition~\ref{d3}.~\fin

The following theorem is proved by D.~Nelson~\cite{Nels}.

\begin{thm} \label{Nels}
Every formula derived in the intuitionistic arithmetic $\mathsf{HA}$ is realizable.
\end{thm}

\begin{prop} \label{p15} 
The formula $Q$ is realizable.
\end{prop}

\proof
The formulas $A_1$-$A_{25}$ are the axioms of the predicate version of Robinson arithmetic $\mathsf{IRA}^*$.
All of them are deducible in the theory $\mathsf{HA}^*$, the predicate version of the theory $\mathsf{HA}$.
The deducibility in the theory $\mathsf{HA}^*$ of the formulas $A_{26}$ and $A_{28}$ follows from~(\ref{f33}) and~(\ref{f21}).
The deducibility of $A_{27}$ in this theory follows from Proposition~\ref{p13} and the connection between the theories $\mathsf{HA}^f$, 
$\mathsf{HA}$, and $\mathsf{HA}^*$ discussed above.
Now the realizability of the formulas $A_1$-$A_{28}$ follows from Theorem~\ref{Nels}.~\fin

From the point of view of {\it constructive semantics}, a statement is true if and only if it is realizable.

An arithmetical formula $\Phi$ is called almost negative if it does not contain $\lor$ and contains $\exists$ only in a combination with an atomic formula.
ECT denotes the scheme 
\begin{equation} \label{f41}
\forall x\,(\Psi(x)\to\exists y\,\Phi(x,y))\to\exists e\forall x\,(\Psi(x)\to\exists y\,(\{e\}(x)=y\land\Phi(x,y))),
\end{equation}
where $\Psi$ is an almost negative formula.
This scheme is called {\it extended Church thesis}.
Extended Church thesis is sound with respect to the constructive semantics in the sense that any arithmetical formula of the form~(\ref{f41}), where $\Psi$ 
is an almost  negative formula, is realizable.
Extended Church thesis is not valid from the point of view of classical semantics: one can find an arithmetical formula $\Phi$ and an almost 
negative arithmetical formula $\Psi$ such that the statement~(\ref{f41}) is false.

Another important law of constructive logic is {\it the principle of constructive selection}, proposed by A.~A.~Markov~\cite{Mark2} and
now called {\it the Markov principle}.
It is expressed by the following scheme M:

\begin{equation} \label{f42}
\forall x\,(\Phi(x)\lor\neg\Phi(x))\,\&\,\neg\neg\exists x\,\Phi(x)\to\exists x\,\Phi(x).
\end{equation}
The Markov principle is sound with respect to the semantics of realizability: any arithmetical formula of the form~(\ref{f42}) is realizable.

A detailed discussion of the extended Church thesis and the Markov principle in the context of constructive arithmetic can be found in the 
monograph~\cite{Drag1}.

In the theory $\mathsf{HA}+\rm M$, the principle ECT is equivalent to the following  scheme nCT:
\begin{equation} \label{f429}
\forall x\,(\neg\Psi(x)\to\exists y\,\Phi(x,y))\to\exists e\forall x\,(\neg\Psi(x)\to\exists y\,(G(e,x,y)\,\&\,\Phi(x,y))),
\end{equation}
where $\Psi(x)$ and $\Phi(x,y)$ are arbitrary arithmetical formulas.

The theory $\mathsf{HA}+\rm M+\rm{ECT}$ (and the equivalent theory $\mathsf{HA}+\rm M+\rm{nCT}$) is usually called {\it Markov arithmetic};
let's denote it~$\mathsf{MA}$.

\section{Constructive predicate logic} \label{sklp}

\subsection{Predicate formulas and schemes} \label{sspf}

{\it The language of predicate logic} $L$ is an elementary language whose signature consists of an infinite set of predicate variables $P^n_i$
($i,n\in\mathbb N$), where $P^n_i$ is called a $n$-ary predicate variable and the number $n$ is called {\it arity} of the predicate
variable~$P^n_i$.
Along with the language $L$, we will consider the elementary language $L^*$, the signature of which is obtained by adding to the signature of the 
language $L$ individual constants $0,1,2,\ldots$ for all natural numbers.
The formulas of the languages $L$ and $L^*$ will be called respectively {\it predicate $L$- and $L^*$-formulas}.

Above, we discussed the principles ECT and nCT calling them schemes and understanding this as a kind of a general form of formulas.
A strict concept of a scheme is introduced by M.~M.~Kipnis~\cite{Kipnis}.
{\it A scheme} over the arithmetical language is a formula of an elementary language whose signature is the union of the signatures of the 
arithmetical language $Ar$ and the language of predicate logic~$L$.
Actually, the scheme is a predicate formula containing predicate constants.
More formally, {\it the scheme language} $LS$ is an elementary language whose signature consists of the predicate symbols $Z,E,S,A,M$ of the 
arithmetical language $Ar$ and the predicate variables $P^n_i$ of the language of predicate logic~$L$.
The inductive definition of the scheme is as follows.

\begin{defn} \label{d45} \rm
1) If $\Phi$ is an atomic $Ar$-formula, the $\Phi$ is a scheme.

2) If $P^n_i$ is a predicate variable, $v_1,\ldots,v_n$ are individual variables, then $P^n_i(v_1,\ldots,v_n)$ is a scheme.

3) If $\Phi$ and $\Psi$ are schemes, then $(\Phi\,\&\,\Psi)$, $(\Phi\lor\Psi)$, $(\Phi\to\Psi)$ are schemes.

4) If $\Phi$ is a scheme, then $\neg\Phi$ is a scheme.

5) If $\Phi$ is a scheme, $v$ is an individual variable, then $\exists v\,\Phi$ and $\forall v\,\Phi$ are schemes.
\end{defn}

Thus arithmetical $Ar$-formulas and predicate $L$-formulas are special cases of schemes.

For technical purposes, we will introduce the language $LS^*$ as an extension of the scheme language $LS$ by adding the individual constants 
$0,1,2,\ldots$ for all natural numbers.
The definition of the $LS^*$-scheme is obtained by replacing in Definition~\ref{d45} points 1) and 2) with the following:

1) If $\Phi$ is an atomic $Ar^*$ formula, then $\Phi$ is a scheme.

2) If $P^n_i$ is a predicate variable, $t_1,\ldots,t_n$ are terms (i.e. individual variables or constants), then $P^n_i(t_1,\ldots,t_n)$ is a
scheme.

If the notation $\Phi(P_1,\ldots,P_n,x_1,\ldots,x_m)$ is used for the scheme $\Phi$, then this means that $\Phi$ does not contain predicate
variables other than $P_1,\ldots,P_n$ and free individual variables other than $x_1,\ldots,x_m$.

Let $\Psi$ be a formula of some elementary language $\Omega$ (for example, $LA$, $Ar$, the language of predicate logic $L$ or the scheme 
language~$LS$).
We will say that $\Psi$ {\it is free} for an $m$-ary predicate variable $P$ in the scheme $\Phi$ if the following
conditions are fulfilled:
1) for any atomic subformula $P(v_1,\ldots,v_m)$ of the scheme $\Phi$, no free occurence of the variable $x_j$ ($j=1,\ldots,m$) in $\Psi$ is within the scope of the quantifier $\exists v_j$ or $\forall v_j$;
2) no atomic subformula $P(v_1,\dots,v_n)$ of the scheme $\Phi$ is not in the scope of the quantifier $\forall v$ or $\exists v$, where $v$ is
a free variable of $\Psi$ other than $x_1,\dots,x_n$.

We will say that a list of formulas $\Psi_1,\dots,\Psi_n$ of an elementary language $\Omega$ {\it is acceptable for substitution} in the scheme 
$\Phi(P_1,\dots,P_n)$ if everyone of the formulas $\Psi_i$ ($i=1,\dots,n$) is free for $P_i$ in $\Phi(P_1,\dots,P_n)$.
By $\Phi(\Psi_1,\dots,\Psi_n)$ denote the formula of the language $\Omega$ obtained by replacing in $\Phi(P_1,\dots,P_n)$ each
atomic subformula of the form $P_i(v_1,\dots,v_{m_i})$ by the result of simultaneous substitution of the variables $v_1,\ldots,v_{m_i}$ into the 
formula $\Psi_i$ for free occurrences of the variables $x_1,\ldots,x_{m_i}$.
The formula $\Phi(\Psi_1,\dots,\Psi_n)$ will be called {\it a substitutional instance} of the scheme $\Phi(P_1,\dots,P_n)$.

Note that if formulas $\Psi_1,\dots,\Psi_m$ are free for predicate variables $P_1,\dots,P_m$ in a closed scheme $\Phi(P_1,\dots,P_m)$, then the 
parameters of $\Phi(\Psi_1,\dots,\Psi_m)$ are exactly \lq\lq superfluous\rq\rq\, parameters of $\Psi_1,\dots,\Psi_m$, i.e. free 
variables of each of the formulas $\Psi_i$ other than $x_1,\dots,x_{m_i}$, where $m_i$ is arity of the predicate variable~$P_i$.

Now we see that the Markov principle is expressed by the following predicate formula~$M$:
$$\forall x\,(P(x)\lor\neg P(x))\,\&\,\neg\neg\exists x\,P(x)\to\exists x\,P(x)$$
in the sense that every formula of the form~(\ref{f42}) is an arithmetical substitutional instance of the formula $M$.
The principle nCT discussed above can be expressed by the following scheme $nCT$ in the strict sense of the word:
$$\forall x\,(\neg P(x)\to\exists y\,Q(x,y))\to\exists z\forall x\,(\neg P(x)\to\exists y\,(G(z,x,y)\,\&\,Q(x,y))),$$
where $P$ and $Q$ are predicate variables, $G(z,x,y)$ is the arithmetical $\Sigma$-formula expressing the predicate $\{z\}(x)=y$.
Indeed, every formula of the form~(\ref{f429}) is an arithmetical substitutional instance of the scheme $nCT$.

\subsection{Absolute realizability} \label{star}

The following notion is introduced by the author in the paper~\cite{Plisko16}.
A motivation is also given there.

\begin{defn} \label{d48} \rm
A $k$-ary {\it generalized predicate} is an arbitrary function (in the set-theoretical sense) of the type
${\mathbb N}^k\to 2^{\mathbb N}$.
\end{defn}

{\it An interpretation} is a partial mapping $f$ which takes some $n$-ary generalized predicate $f(P^n_i)$ to each predicate variable $P^n_i$ from 
its domain.
We will say that $f$ is an interpretation of the scheme $\Phi$ if $f$ is defined on all predicate variables in~$\Phi$.
For technical purposes, we will introduce the predicate constant $\bot$ into the language, which will be considered an atomic formula.

\begin{defn} \label{d28} \rm
Let $f$ be an interpretation of a closed $LS^*$-scheme~$\Phi$.
The relation $e\rf\Phi$, where $e\in\mathbb N$, is defined by induction on the number of logical symbols in~$\Phi$.

\begin{itemize}
\item
It is not true that $e\rf\bot$, whatever $e\in\mathbb N$ is.
\item
$\rf\Phi\rightleftharpoons e\rk\Phi$ if $\Phi$ is an atomic $Ar^*$-formula.
\item
$e\rf P(a_1,\ldots,a_n)\rightleftharpoons e\in f(P)(a_1,\ldots,a_n)$ if $P$ is an $n$-ary predicate variable, $a_1,\ldots,a_n\in\mathbb N$.
\item
$e\rf(\Phi\,\&\,\Psi)\rightleftharpoons(e)_0\rf\Phi$ and $(e)_1\rf\Psi$.
\item
$e\rf(\Phi\lor\Psi)\rightleftharpoons(e)_0=0$ and $(e)_1\rf\Phi$ or $(e)_0=1$ and $(e)_1\rf\Psi$.
\item
$e\rf(\Phi\to\Psi)\rightleftharpoons\forall a\,(a\rf\Phi\Rightarrow!\{e\}(a)$ and $\{e\}(a)\rf\Psi)$.
\item
$e\rf\neg\Phi\rightleftharpoons e\rf(\Phi\to\bot)$.
\item
$e\rf\exists x\,\Phi(x)\rightleftharpoons(e)_1\rf\Phi((e)_0)$.
\item
$e\rf\forall x\,\Phi(x)\rightleftharpoons\forall a\,(!\{e\}(a)$ and $\{e\}(a)\rf\Phi(a))$.
\end{itemize}
\end{defn}

If $e\rf\Phi$ holds, then we will say that {\it the natural number $e$ realizes the scheme $\Phi$ in the interpretation $f$} or is
{an \it $f$-realization of the scheme $\Phi$}.
We will say that a closed scheme $\Phi$ is {\it realizable in the interpretation $f$}, or is $f$-{\it realizable} if there is
a natural number $e$ such that $e\rf\Phi$ holds.

It follows from Definitions~\ref{d3} and~\ref{d28} that if $\Phi$ is a closed $Ar^*$ formula, then 
$\forall e\,[e\rf\Phi\Leftrightarrow e\rk\Phi]$ for any interpretation~$f$.

There is the following theorem about the correctness of the intuitionistic predicate calculus $\mathsf{IQC}$ with respect to $f$-realizability.

\begin{thm} \label{t32} 
If a closed predicate $L^*$-formula $\Phi$ is deducible in $\mathsf{IQC}$, then $\Phi$ is $f$-realizable for any interpretation~$f$.
Moreover, by deducing the predicate $L^*$-formula $\Phi$ in $\mathsf{ICQ}$, one can effectively find a number $e$ such that 
$e\rf\Phi$ for any interpretation~$f$.
\end{thm}

\proof
This theorem is proved in the same way as Nelson's theorem on the correctness of $\mathsf{IQC}$ with respect to recursive
realizability (see~\cite{Nels},~\cite[Theorem~62]{Kleene3}).~\fin

The concepts of a uniformly absolutely realizable and absolutely irrefutable scheme introduced below go back to the semantics of realizability for
predicate formulas based on the interpretation of predicate variables only by predicates defined in the arithmetical language and the semantics
of realizability for this language (see~\cite{Plisko14}).

\begin{defn} \label{d53} \rm
A closed scheme $\Phi$ is called {\it uniformly absolutely realizable} if there exists a natural number $e$ such that
$e\rf\Phi$ for any interpretation~$f$.
\end{defn}

Let $z$ be an individual variable that does not occur in the $LS^*$-scheme~$\Phi$.
By $\Phi[z]$ we denote the scheme obtained by replacing in $\Phi$ each atomic subformula $P^n_i(t_1,\ldots,t_n)$, where
$t_1,\ldots,t_n$ are terms, i.e. individual variables or constants, by $P^{n+1}_i(z,t_1,\ldots,t_n)$.
Obviously, the scheme $\Phi[z]$ contains the free variable $z$, as well as all the free variables of the predicate scheme~$\Phi$.
In particular, if $\Phi$ is a closed scheme, then $z$ is the only parameter of the scheme~$\Phi[z]$.
By $\Phi[t]$ we will denote the result of substituting the term $t$ in $\Phi[z]$ for the free occurrences of~$z$.

If $f$ is an interpretation, $c$ is a natural number, then by $f_c$ we denote the interpretation that takes each predicate variable $P^n_i$
to  the generalized predicate $f_c(P^n_i)$ defined as follows:
$$f_c(P^n_i)(k_1,\ldots,k_n)=f(P^{n+1}_i)(c,k_1,\ldots,k_n).$$

\begin{prop} \label{p50} 
For any interpretation $f$, natural numbers $c,e$, and closed $LS^*$-scheme~$\Phi$,
\begin{equation} \label{f9}
e\rf\Phi[c]\Leftrightarrow e\,\mathsf r^{f_c}\,\Phi.
\end{equation}
\end{prop}

\proof
Induction on the number of logical symbols in the scheme~$\Phi$.
The case when $\Phi$ is $\bot$ is trivial.

If $\Phi$ is an atomic $Ar^*$-formula, the statement is obvious because then $\Phi[c]$ coincides with $\Phi$, thus $e\rf\Phi[c]$ and
$e\,\mathsf r^{f_c}\,\Phi$ mean the same, namely $e\rk\Phi$.

If $\Phi$ is atomic $L^*$-formula $P^n_i(k_1,\ldots,k_n)$, where $k_1,\ldots,k_n$ are constants, then $\Phi[c]$ is the atomic $L^*$-formula
$P^{n+1}_i(c,k_1,\ldots,k_n)$, so we have:
$$e\rf\Phi[c]\Leftrightarrow e\in f(P^{n+1}_i)(c,k_1,\ldots,k_n)\Leftrightarrow e\in f_c(P^n_i)(k_1,\ldots,k_n)\Leftrightarrow 
e\,\mathsf r^{f_c}\,\Phi.$$

We prove that if the condition~(\ref{f9}) is fulfilled when $\Phi$ is an $LS^*$-scheme $\Phi_1$ or $\Phi_2$ for any natural numbers
$c,e$, then this condition is satisfied when $\Phi$ is an $LS^*$-scheme of the form $\Phi_1\circ\Phi_2$, where $\circ$ is $\&$, $\lor$ or $\to$.
Let $\circ$ be $\&$.
Then $\Phi[c]$ is $\Phi_1[c]\,\&\,\Phi_2[c]$, and we have:
$$e\rf\Phi[c]\Leftrightarrow e\rf(\Phi_1[c]\,\&\,\Phi_2[c])\Leftrightarrow(e)_0\rf\Phi_1[c]\mbox{ and }(e)_1\rf\Phi_2[c]\Leftrightarrow$$ 
$$\Leftrightarrow (e)_0\,\mathsf r^{f_c}\,\Phi_1\mbox{ and }(e)_1\,\mathsf r^{f_c}\,\Phi_2\Leftrightarrow e\,\mathsf r^{f_c}\,(\Phi_1\,\&\,\Phi_2)
\Leftrightarrow e\,\mathsf r^{f_c}\,\Phi.$$
The case when $\circ$ is $\lor$ is considered in exactly the same way.

Let $\circ$ be $\to$.
Then $\Phi[c]$ is $\Phi_1[c]\to\Phi_2[c]$, and we have:
$$e\rf\Phi[c]\Leftrightarrow e\rf(\Phi_1[c]\to\Phi_2[c])\Leftrightarrow\forall a\,(a\rf\Phi_1[c]\Rightarrow\{e\}(a)\rf\Phi_2[c])\Leftrightarrow$$ 
$$\Leftrightarrow\forall a\,(a\mathsf r^{f_c}\Phi_1\Rightarrow\{e\}(a)\mathsf r^{f_c}\Phi_2)\Leftrightarrow e\,\mathsf r^{f_c}\,(\Phi_1\to\Phi_2)
\Leftrightarrow e\,\mathsf r^{f_c}\,\Phi.$$

Assume that the statement~(\ref{f9}) holds.
We prove that then it is also fulfilled in the case when the $LS^*$-scheme $\neg\Phi$ is considered in the role of $\Phi$.
We have:
$$e\rf\neg\Phi[c]\Leftrightarrow e\rf(\Phi[c]\to\bot)\Leftrightarrow e\,\mathsf r^{f_c}\,(\Phi\to\bot)\Leftrightarrow 
e\,\mathsf r^{f_c}\,\neg\Phi.$$

Let $\Phi$ be of the form $\exists x\,\Psi(x)$.
Then $\Phi[c]$ is the $LS^*$-scheme $\exists x\,\Psi(x)[c]$.
We prove that the statement~(\ref{f9}) holds if it holds for any $LS^*$-scheme of the form $\Psi(k)$, where $k$ is an arbitrary natural
number.
We have:
$$e\rf\Phi[c]\Leftrightarrow e\rf\exists x\,\Psi(x)[c]\Leftrightarrow (e)_1\rf\Psi((e)_0)[c]\Leftrightarrow$$
$$\Leftrightarrow(e)_1\,\mathsf r^{f_c}\,\Psi((e)_0)\Leftrightarrow e\,\mathsf r^{f_c}\,\exists x\,\Psi(x)\Leftrightarrow 
e\,\mathsf r^{f_c}\,\Phi.$$
The case when $\Phi$ is of the form $\forall x\,\Psi(x)$ is considered similarly.~\fin

\begin{defn} \label{d26} \rm
A closed $LS^*$-scheme $\Phi$ is called {\it absolutely irrefutable} if the $LS^*$-scheme $\forall z\,\Phi[z]$ is $f$-realizable for any
interpretation~$f$.
\end{defn}

Let's agree on the following designations.
Let $\varphi(x_1,\ldots,x_n)$ be an expression specifying an $n$-ary partial recursive function.
Then we can find a G\"odel number of this function denoted by $\Lambda x_1,\ldots,x_n.\varphi(x_1,\ldots,x_n)$.
If the expression $\varphi(x_1,\ldots,x_n,y_1,\ldots,y_m)$ defines an $(m+n)$-place partial recursive function $\varphi$, then the expression
$\Lambda x_1,\ldots,x_n.\varphi(x_1,\ldots,x_n,y_1,\ldots,y_m)$ denotes an $m$-place recursive function.

\begin{thm} \label{t7}
A closed $LS^*$-scheme is absolutely uniformly realizable if and only if it is absolutely irrefutable.
\end{thm}

\proof
Let $\Phi$ the a closed absolutely uniformly realizable $LS^*$-scheme, i.e. there exists a natural number $e$ such that
$e\rf\Phi$ holds for any interpretation~$f$.
We prove that $\Phi$ is absolutely irrefutable, i.e. that the $LS^*$-scheme $\forall z\,\Phi[z]$ is $f$-realizable for any interpretation~$f$.
Let $d=\Lambda z.e$.
We prove that $d\rf\forall z\,\Phi[z]$ for any interpretation~$f$.
It is required to prove that $!\{d\}(c)$ for any $c$ (this condition is obviously fulfilled) and $\{d\}(c)\rf\Phi[c]$, i.e.
$e\rf\Phi[c]$.
But this is obvious since by Proposition~\ref{p50}, it is equivalent to $e\,\mathbf r^{f_c}\,\Phi$, and the latter 
holds  by the condition.

Prove the converse: if a closed $LS^*$-scheme $\Phi$ is absolutely irrefutable, then $\Phi$ is absolutely uniformly realizable.
Let $\Phi$ be an absolutely irrefutable closed $LS^*$-scheme.
Suppose however that it is not absolutely uniformly realizable.
The latter means that for every natural number $k$ there is an interpretation $f^k$ such that $k\,\mathbf r^{f^k}\,\Phi$ does not hold.

The interpretation $f$ of the $LS^*$-scheme $\forall z\,\Phi[z]$ is defined as follows: if the predicate variable $P^n_i$ occurs in $\Phi$, then 
for any $k,k_1,\ldots,k_n$ we put
\begin{equation} \label{f61}
f(P^{n+1}_i)(k,k_1,\ldots,k_n)=
\begin{cases}
f^{\{k\}(k)}(P^n_i)(k_1,\ldots,k_n)&\mbox{if }!\{k\}(k);\\
\emptyset&\mbox{otherwise}.
\end{cases}
\end{equation}
We prove that the $LS^*$-scheme $\forall z\,\Phi[z]$ is not $f$-realizable.
This would mean that $\Phi$ is not absolutely irrefutable contrary to the condition.
Assume that $e\rf\forall z\,\Phi[z]$ for some~$e$.
Then $!\{e\}(k)$ for every~$k$.
In~particular, $!\{e\}(e)$, and at the same time $\{e\}(e)\rf\Phi[e]$.
By Proposition~\ref{p50}, we have
\begin{equation} \label{f17}
\{e\}(e)\,\mathsf r^{f_e}\Phi.
\end{equation}
Note that for each predicate variable $P^n_i$ in $\Phi$ and any $k_1,\ldots,k_n$, the condition
$$f_e(P^n_i)(k_1,\ldots,k_n)=f(P^{n+1}_i)(e,k_1,\ldots,k_n)=f^{\{e\}(e)}(P^n_i)(k_1,\ldots,k_n)$$
is fulfilled.
This means that the interpretations $f_e$ and $f^{\{e\}(e)}$ coincide on all predicate variables in~$\Phi$.
Then~(\ref{f17}) implies $\{e\}(e)\,\mathsf r^{f^{\{e\}(e)}}\Phi$, i.e. $k\,\mathsf r^{f^k}\,\Phi$ for $k=\{e\}(e)$
contrary to the main property of the interpretation~$f^k$.~\fin

Thus the concepts of an absolutely uniformly realizable and absolutely irrefutable $LS^*$-scheme coincide.

\begin{defn} \label{d64} \rm
A closed $LS^*$-scheme is called {\it absolutely realizable} if it is absolutely uniformly realizable or, equivalently, if
it is absolutely irrefutable.
\end{defn}

If $\Phi(x_1,\ldots,x_n)$ is a scheme, then by $\tilde\forall\Phi$ we will denote {\it the universal closure} of the scheme $\Phi$, namely
the scheme $\forall x_1,\ldots,x_n\,\Phi$.

\begin{prop} \label{p54} 
If a closed scheme $\Phi(P_1,\ldots,P_n)$ is absolutely realizable, then the scheme $\tilde\forall\,\Phi(\Psi_1,\dots,\Psi_n)$ 
is absolutely realizable for any list of schemes $\Psi_1,\dots,\Psi_n$ acceptable for substitution in $\Phi(P_1,\ldots,P_n)$.
\end{prop}

\proof
Let $\Phi(P_1,\ldots,P_n)$ be an absolutely realizable closed scheme.
This means that there is a number $e$ such that $e\rf\Phi$ for any interpretation~$f$.
Assume that schemes $\Psi_1,\ldots,\Psi_n$ are given and $y_1,\ldots,y_s$ are the free variables of the scheme $\Phi(\Psi_1,\ldots,\Psi_n)$, so 
this scheme can be denoted as $\Phi(\Psi_1,\ldots,\Psi_n,y_1,\ldots,y_s)$.
Note also that each of the schemes $\Psi_i$ ($i=1,\ldots,n$) contains no free variables except $y_1,\ldots,y_s$ and $x_1,\ldots,x_{m_i}$, where
$m_i$ is arity of the predicate variable $P_i$, so for $\Psi_i$ it is appropriate to use the notation $\Psi_i(y_1,\ldots,y_s,x_1,\ldots,x_{m_i})$.
Let's prove that the number $\Lambda y_1.\ldots\Lambda y_m.e$ is an $f$-realization of the scheme $\tilde\forall\,\Phi(\Psi_1,\ldots,\Psi_n)$ for 
any interpretation~$f$.
This means that for any numbers $k_1,\ldots,k_m$,
\begin{equation} \label{f65}
e\rf\Phi(\Psi_1,\ldots,\Psi_n,k_1,\ldots,k_m).
\end{equation}
The interpretation $g$ is defined as follows: for each $i=1,\ldots,n$ and any natural $a,\ell_1,\ldots,\ell_{m_i}$ we put
$$a\in g(P_i)(\ell_1,\ldots,\ell_{m_i})\Leftrightarrow a\rf\Psi_i(k_1,\ldots,k_s,\ell_1,\ldots,\ell_{m_i}).$$
It is not difficult to make sure that for any number $b$,
$$b\,\mathsf r^g\,\Phi(P_1,\ldots,P_n)\Leftrightarrow b\rf\Phi(\Psi_1,\ldots,\Psi_n,k_1,\ldots,k_m).$$
Since by the condition, $e\,\mathsf r^g\,\Phi(P_1,\ldots,P_n)$, we have~(\ref{f65}), as was to be proved.~~\fin

Proposition~\ref{p54} means that the universal closure of a substitutional instance of an absolutely realizable predicate formula is absolutely
realizable.

\begin{thm} \label{t18} 
The predicate formula
$$\forall x\,(P(x)\lor\neg P(x))\,\&\,\neg\neg\exists x\,P(x)\to\exists x\,P(x)\eqno{(M)}$$
is absolutely realizable.
\end{thm}

\proof
Let $e=\Lambda x.2^{l(x)}\cdot3^{r(x)}$, where
$$l(x)=\mu y[(\{(x)_0\}(y))_0=0],\;r(x)=(\{(x)_0\}(l(x)))_1.$$
We prove that $e$ realizes $(M)$ in any interpretation~$f$.
The latter means that for any natural $a$, if
\begin{equation} \label{f19}
a\rf\forall x\,(P(x)\lor\neg P(x))\,\&\,\neg\neg\exists x\,P(x),
\end{equation}
then $2^{l(a)}\cdot3^{r(a)}\rf\exists x\,P(x)$, i.e.
\begin{equation}\label{f196}
r(a)\rf P(l(a)).
\end{equation}
So let $a$ be such that~(\ref{f19}) holds.
Then
\begin{equation} \label{f198}
(a)_0\rf\forall x\,(P(x)\lor\neg P(x))
\end{equation}
and $(a)_1\rf\neg\neg\exists x\,P(x)$.
By the principle of constructive selection, it follows that there exists a number $b$ such that
$(\{(a)_0\}(b))_0=0$.
Obviously, the least such number $l(a)$ can be found effecively.
It follows from the condition~(\ref{f198}) that $(\{(a)_0\}(l(a)))_1\rf P(l(a))$, i.e.~(\ref{f196}) holds, as was to be proved.~\fin

\section{Scheme theorem} \label{stos}

\subsection{Predicate $Ar$-formulas}

In the previous section, we defined the concept of an absolutely realizable scheme and, in particular, an absolutely realizable predicate formula.
We see that actually a scheme is a predicate formula containing predicate symbols with a fixed rigid interpretation.
The {\it scheme theorem} proven in~\cite{Plisko14} in the context of recursive realizability allows us to get rid of the rigid interpretation of predicate constants and treat them as predicate variables.
Our immediate goal is to prove a variant of the scheme theorem in the context of absolute realizability.
The following concepts reflect the idea of simultaneously considering some predicate symbols both as predicate constants and as predicate
variables.

By $Z$ denote the predicate variable $P_0^1$, by $E$ and $S$ denote respectively the predicate variables $P_0^2$ and
$P_1^2$, by $A$ and $M$ denote the predicate variables $P_0^3$ and $P_1^3$ respectively.
In these notations, the languages $Ar$ and $Ar^*$ turn out to be fragments of the languages $L$ and $L^*$ respectively.
Formulas of the language $L$ or $L^*$ that do not contain predicate variables other than $Z,E,S,A,M$ we will call respectively predicate $Ar$- or $Ar^*$-formulas.

Thus each formula of the language $Ar$ ($Ar^*$) can be considered as a predicate formula of the language $L$ ($L^*$) if the predicate
symbols $Z,E,S,A,M$ are treated as predicate variables.
Since a predicate $Ar$- or $Ar^*$-formula $\Phi$ does not contain predicate variables other than $Z,E,S,A,M$, we will denote it
$\Phi(Z,E,S,A,M)$.
An interpretation $f$ is called {\it a model} of a closed predicate $Ar$-formula $\Phi$ if there exists a natural number $e$ such that
$e\rf\Phi$.

In Section~\ref{ska}, the arithmetical formula $Q$ was defined as the conjunction of the formulas $A_1$-$A_{28}$.
As noted above, $A_1$-$A_{25}$ set the axiomatics of the predicate variant of Robinson's intuitionistic arithmetic $\mathsf{IRA}$.
Therefore, if an $LA$-formula is deduced in $\mathsf{IRA}$, then its predicate form is deduced in $\mathsf{IQC}$ from
the hypothesis~$Q$.
This leads us to the following statement.

\begin{prop} \label{p55} 
If an interpretation $f$ is a model of the predicate $Ar$-formula $Q$, then for any closed $LA$-formula $\Phi$ deducible in 
$\mathsf{IRA}$, it is possible to effectively construct an $f$-realization of the predicate form of~$\Phi$.
\end{prop}

\proof
Suppose an interpretation $f$ is a model of the predicate $Ar$-formula $Q$, i.e  there is a number $a$ such that $a\rf Q$.
Suppose $\mathsf{IRA}\vdash\Phi$.
Then $Q\vdash\Phi^*$, where $\Phi^*$ is the predicate form of $\Phi$ and $\vdash$ means deducibility in $\mathsf{ICQ}$.
By the deduction theorem, it is possible to construct the derivation of $Q\to\Phi^*$ in $\mathsf{ICQ}$, and by 
Theorem~\ref{t32}, one can find a number $b$ such that $b\rf (Q\to\Phi^*)$, and then $\{a\}(b)\rf\Phi^*$.~\fin

If $\vec  x$ is the list of variables $x_1,\ldots,x_n$, $\vec y$ is the list $y_1,\ldots,y_n$, then $E(\vec x,\vec y)$ will  denote the formula $E(x_1,y_1)\,\&\ldots\&\,E(x_n,y_n)$.

\begin{prop} \label{p57}
If an interpretation $f$ is a model of the predicate $Ar$-formula $Q$, then for any predicate $Ar^*$-formula $\Phi(\vec z,\vec x)$, the formula
$$\forall\vec z\forall\vec x,\vec y\,(E(\vec x,\vec y)\to(\Phi(\vec z,\vec x)\to\Phi(\vec z,\vec y)))$$
is $f$-realizable.
\end{prop}

\proof
The statement follows from Proposition~\ref{p55}, the fact that the formula
$$\forall\vec z\forall\vec  x,\vec y\,(x_1=y_1\,\&\ldots\&\,x_n=y_n\to(\Psi(\vec z,\vec x)\to\Psi(\vec z,\vec y))$$
is deducible in $\mathsf{IRA}$ for any $LA$-formula $\Psi(\vec x)$, and a simple remark that
every predicate $Ar$-formula $\Phi(\vec x)$ is a predicate form of some $LA$-formula $\Psi(\vec x)$.~\fin

\begin{prop} \label{p56}
If an interpretation $f$ is a model of the predicate $Ar$-formula $Q$, then for any predicate $Ar^*$-formula $\Phi(\vec x)$ there is an $(n+2)$-place
partial recursive function $\gamma_{\Phi}$ such that for any natural $d,e,\vec k,\vec\ell$, if $d\rf E(\vec k,\vec\ell)$ and $e\rf\Phi(\vec k)$, then
$\gamma_\Phi(d,e,k,\vec\ell)\rf\Phi(\ell)$.
\end{prop}

\proof
This is a simple consequence of Proposition~\ref{p57}.~\fin

\subsection{Standard elements} \label{sse}
For a natural $n$, the predicate $Ar$-formula $[n](x)$ is defined inductively as follows:
\begin{itemize}
\item
$[0](x)$ is $Z(x)$;
\item
$[n+1](x)$ is $\exists x_n\,([n](x_n)\,\&\,S(x_n,x))$.
\end{itemize}

Note that $[n](x)$ is the predicate form of the $LA$-formula $x=n$.

\begin{prop} \label{p118} 
If an interpretation $f$ is a model of the predicate $Ar$-formula $Q$, then for every $f$-realization of $Q$ and every natural number 
$n$, it is possible to effectively construct $f$-realizations of the following predicate $Ar$-formulas:
\begin{equation} \label{f55}
\forall x,y\,([n](x)\,\&\,[n](y)\to E(x,y));
\end{equation}
\begin{equation} \label{f58}
\forall x\,([n](x)\lor\neg[n](x));
\end{equation}
\begin{equation} \label{f585}
\exists x\,[n](x);
\end{equation}
\begin{equation} \label{f589}
\forall x\,([n](x)\to\neg[m](x))
\end{equation}
if $m\not=n$.
\end{prop}

\proof
Every one of the formulas~(\ref{f55})-(\ref{f589}) is the predicate form of some $LA$-formula derived in $\mathsf{IRA}$,
thus the statement follows from Proposition~\ref{p55}.~\fin

\begin{prop} \label{p61}
If an interpretation $f$ is a model of the predicate $Ar$-formula $Q$, then for any natural numbers $n$, $m$, and $\ell$, if the predicate 
$Ar^*$-formulas $[n](\ell)$ and $[m](\ell)$ are $f$-realizable, then $m=n$.
\end{prop}

\proof
Suppose $[n](\ell)$ and $[m](\ell)$ are $f$-realizable.
Assume that $m\not=n$.
Then by Proposition~\ref{p118}, the predicate $Ar$-formula~(\ref{f589}) is $f$-realizable, and this implies $f$-realizability of the 
predicate $Ar^*$-formula $\neg[m](\ell)$ contrary to $f$-realizability of $[m](\ell)$.~\fin

\begin{prop} \label{p59} 
If an interpretation $f$ is a model of the predicate $Ar$-formula $Q$, then for any $LA$-formula $\Phi(x_1,\ldots,x_n)$ and any natural numbers 
$k_1,\ldots,k_n$, if $\mathsf{IRA}\vdash\Phi(k_1,\ldots,k_n)$, 
then for every $f$-realization of $Q$ one can effectively construct an $f$-realization of the formula
\begin{equation} \label{f62}
\forall x_1,\ldots,x_n\,([k_1](x_1)\,\&\ldots\&\,[k_n](x_n)\to\Psi(x_1,\ldots,x_n)),
\end{equation}
where $\Psi$~is the predicate form of~$\Phi$.
\end{prop}

\proof
This follows from Proposition~\ref{p55} and the fact that $\mathsf{IRA}\vdash\Phi(k_1,\ldots,k_n)$ implies
$\mathsf{IRA}\vdash\forall x_1,\ldots,x_n\,(x_1=k_1\,\&\ldots\&\,x_n=k_n\to\Phi(x_1,\ldots,x_n)).$~\fin

\begin{prop} \label{p62} 
If an interpretation $f$ is a model of the predicate $Ar$-formula $Q$, then for every $f$-realization of $Q$ and every natural 
number $n$, it is possible to effectively construct natural numbers $\tilde n$ and $e$ such that $e$ is an $f$-realization of the predicate 
$Ar^*$-formula $[n](\tilde n)$.
\end{prop}

\proof
By Proposition~\ref{p118}, for any $f$-realization of $Q$ and any natural number $n$, it is possible to effectively construct 
an $f$-realization $a$ of~(\ref{f585}).
Then $(a)_1\rf[n]((a)_0)$, so one can put $\tilde n=(a)_0$, $e=(a)_1$.~\fin

Let's fix an algorithm that takes each $f$-realization of $Q$ and each number $n$ to the number~$\tilde n$.

\begin{prop} \label{p64} 
Suppose an interpretation $f$ is a model of the formula $Q$; then for any arithmetical formula $\Phi(x_1,\ldots,x_n)$ and natural 
numbers $k_1,\ldots,k_n$, it follows from  $\mathsf{IRA}\vdash\Phi(k_1,\ldots,k_n)$ that for every $f$-realization of $Q$, it is 
possible to effectively construct an $f$-realization of the predicate $Ar^*$-formula 
$\Psi(\tilde k_1,\ldots,\tilde k_n)$, where $\Psi$ is the predicate form of~$\Phi$.
\end{prop}

\proof
Suppose $\mathsf{IRA}\vdash\Phi(k_1,\ldots,k_n)$.
By Proposition~\ref{p59}, for any $f$-realization of the predicate $Ar$-formula $Q$, it is possible to effectively construct an $f$-realization of the 
predicate $Ar$-formula~(\ref{f62}).
After that, one can effectively find an $f$-realization of the predicate $Ar^*$-formula
$$[k_1](\tilde k_1)\,\&\ldots\&\,[k_n](\tilde k_1)\to\Psi(\tilde k_1,\ldots,\tilde k_n).$$
By Proposition~\ref{p62}, an $f$-realization of the premise of this formula can be found, and then an $f$-realization of the conclusion, i.e.
the predicate $Ar^*$-formula $\Psi(\tilde k_1,\ldots,\tilde k_n)$, is constructed.~\fin

\begin{prop} \label{p66} 
Suppose an interpretation $f$ is a model of $Q$.
Then for any arithmetical $\Sigma$-formula $\Phi(x_1,\ldots,x_n)$ and natural $k_1,\ldots,k_n$, if the $Ar^*$-formula 
$\Phi(k_1,\ldots,k_n)$ is true, then the predicate $Ar^*$-formula $\Phi(\tilde k_1,\ldots,\tilde k_n)$ is $f$-realizable and its 
$f$-realization can be found effectively.
\end{prop}

\proof
It follows from the results on the numerical expressibility~\cite[\S~41]{Kleene3} that every true $\Sigma$-statement is deducible in $\mathsf{IRA}$, thus the proposition being proved follows from Proposition~\ref{p64}.\fin

Recall that the $Ar$-formula $G(x,y,z)$ expresses the predicate $\{x\}(y)=z$.

\begin{prop} \label{lemma7}
Suppose an interpretation $f$ is a model of the predicate $Ar$-formula $Q$ and the function $\{e\}$ is total; then for any natural numbers $n$ and 
$k$, the predicate $Ar^*$-formula $G(\tilde e,\tilde n,\tilde k)$ is $f$-realizable if and only if $\{e\}(n)=k$.
\end{prop}

\proof
Suppose $\{e\}$ is a total function, and $\{e\}(n)=k$.
Then $Ar^*$-formula $G(e,n,k)$ is true and by Proposition~\ref{p66}, the predicate $Ar^*$-formula $G(\tilde e,\tilde n,\tilde k)$
is $f$-realizable.
Conversely, suppose the predicate $Ar^*$-formula $G(\tilde e,\tilde n,\tilde k)$ is $f$-realizable.
Since $\{e\}$ is a total function, there exists a natural number $\ell$ such that $\{e\}(n)=\ell$.
Then $G(e,n,\ell)$ is true.
Reasoning as above, we conclude that the predicate $Ar^*$-formula $G(\tilde e,\tilde n,\tilde\ell)$ is $f$-realizable.
Since the predicate $Ar$-formula $A_{26}$ is $f$-realizable, it follows that the predicate $Ar^*$-formula $E(\tilde\ell,\tilde k)$ is 
$f$-realizable.
The predicate $Ar^*$-formula $[k](\tilde k)$ is $f$-realizable.
It follows from Proposition~\ref{p57} that the predicate $Ar^*$-formula $[k](\tilde\ell)$ is $f$-realizable.
On the other hand, the predicate $Ar^*$-formula $[\ell](\tilde\ell)$ is $f$-realizable.
Then by Proposition~\ref{p61}, $k=\ell$.
Thus we proved that $\{e\}(n)=k$.~\fin

Recall that $H(x_1,x_2)$ is the $\Sigma$-formula $\exists z\,(Z(z)\,\&\,G(x_1,x_2,z)).$

\begin{prop} \label{p122}
Suppose an interpretation $f$ is a model of the predicate $Ar$-formula $Q$ and the function $\{e\}$ is total; then for any natural number $n$, the 
predicate $Ar^*$-formula $H(\tilde e,\tilde n)$ is $f$-realizable if and only if $\{e\}(n)=0$.
\end{prop}

\proof
Suppose $\{e\}(n)=0$.
By Proposition~\ref{lemma7}, the predicate $Ar^*$-formula $G(\tilde e,\tilde n,\tilde 0)$ is $f$-realizable.
On the other hand, the predicate $Ar^*$-formula $Z(\tilde 0)$ is also $f$-realizable.
Then the predicate $Ar^*$-formula $Z(\tilde0)\,\&\,G(\tilde e,\tilde n,\tilde0)$ is $f$-realizable.
It follows that the predicate $Ar^*$-formula $\exists z\,(Z(z)\,\&\,G(\tilde e,\tilde n,z))$, i.e. $H(\tilde e,\tilde n)$, is $f$-realizable.

Conversely, suppose the predicate $Ar^*$-formula $H(\tilde e,\tilde n)$ is $f$-realizable.
Then there exists a natural number $k$ such that the predicate $Ar^*$-formulas $Z(k)$, i.e. $[0](k)$, and $G(\tilde e,\tilde n,k)$ are 
$f$-realizable.
Since the predicate $Ar^*$-formula $[0](\tilde 0)$ also is $f$-realizable, then it follows from Proposition~\ref{p118} that the predicate
$Ar^*$-formula $E(k,\tilde0)$ is $f$-realizable.
Then it follows from Proposition~\ref{p57} that $G(\tilde e,\tilde n,\tilde0)$ is $f$-realizable.
By Proposition~\ref{lemma7}, $\{e\}(n)=0$.~\fin

Assume that an interpretation $f$ is a model of the predicate $Ar^*$-formula $Q$.
A natural number $a$ will be called $f$-{\it standard} if there exists a natural number $n$ such that the predicate $Ar^*$-formula $[n](a)$  is  $f$-realizable.
Note that the numbers $\tilde 0,\tilde 1,\tilde 2,\ldots$ are $f$-standard.

\begin{prop} \label{p154}
If a natural number $a$ is $f$-standard, $b$ is a natural number such that the predicate $Ar^*$-formula $S(a,b)$ is $f$-realizable,
then the number $b$ is $f$-standard.
\end{prop}

\proof
Suppose the predicate $Ar^*$-formulas $[n](a)$ and $S(a,b)$ are $f$-realizable.
Then the predicate $Ar^*$-formula $[n](a)\,\&\,S(a,b)$ is $f$-realizable, therefore, the predicate $Ar^*$-formula
$\exists x_n\,([n](x_n)\,\&\,S(x_n,b))$, i.e. $[n+1](b)$, is $f$-realizable.
This means that the number $b$ is $f$-standard.~\fin

In accordance with the notation introduced above, $x\le y$ is the predicate $Ar$-formula $\exists z\,A(z,x,y)$.

\begin{prop} \label{lemma8} 
Whatever natural numbers $n$ and $b$ are, if $b$ is not $f$-standard, then the predicate $Ar^*$-formula $\tilde n\le b$ is $f$-realizable.
\end{prop}

\proof
Induction on~$n$.
Assume $n=0$.
It follows from the $f$-realizability of the predicate $Ar$-formula $A_5$ that the predicate $Ar^*$-formula $Z(\tilde0)\to A(b,\tilde0,b)$ is 
$f$-realizable.
Since the predicate $Ar^*$-formula $Z(\tilde0)$ is $f$-realizable, it follows that the predicate $Ar^*$-formula $A(b,\tilde0,b)$ is $f$-realizable.
Therefore, the predicate $Ar^*$-formula $\exists z\,A(z,\tilde0,b)$, i.e. $\tilde0\le b$, is $f$-realizable.

Now suppose that for any natural $b$ that is not $f$-standard, the predicate $Ar^*$-formula $\tilde n\le b$ is $f$-realizable.
Put $n'=n+1$ and prove that for any such number $b$ the predicate $Ar^*$-formula $\widetilde{n'}\le b$ is $f$-realizable.
Since the number $b$ is not $f$-standard, it follows that the predicate $Ar^*$-formula $Z(b)$ is not $f$-realizable.
Since the predicate $Ar$-formula $A_{13}$ is $f$-realizable, it follows that the predicate $Ar^*$-formula $\exists y\,S(y,b)$ is $f$-realizable.
This means the existence of a natural number $a$ such that the predicate $Ar^*$-formula $S(a,b)$ is $f$-realizable.
By Proposition~\ref{p154}, the number $a$ cannot be $f$-standard.
By the induction hipothesis, the predicate $Ar^*$-formula $\exists z\,A(z,\tilde n,a)$ is $f$-realizable.
This means the existence of a natural number $c$ such that $A(c,\tilde n,a)$ is $f$-realizable.
Since $Ar^*$-formula $S(n,n')$ is true, it follows from Proposition~\ref{p66} that the predicate $Ar^*$-formula $S(\tilde n,\widetilde{n'})$
is $f$-realizable.
Since the predicate $Ar$-formula $A_6$ is $f$-realizable, it follows that the predicate $Ar^*$-formula $A(c,\widetilde{n'},b)$ is $f$-realizable.
Therefore, the predicate $Ar^*$-formula $\exists z\,A(z,\widetilde{n'},b)$, i.e. $\widetilde{n'}\le b$, is $f$-realizable.~\fin

\begin{thm} \label{lemma11}
If an interpretation $f$ is a model of the predicate $Ar$-formula $Q$, then every natural number is $f$-standard.
\end{thm}

\proof
Suppose an interpretation $f$ is a model of the predicate $Ar$-formula $Q$.
Then the predicate $Ar$-formula~$A_{28}$ is $f$-realizable.
This implies the existence of a binary recursive function $g$ such that for any $k,\ell$, $g(k,\ell)\rf(B(k,\ell)\lor\neg B(k,\ell))$, i.e., if 
$(g(k,\ell))_0=0$, then the predicate $Ar^*$-formula $B(k,\ell)$ is $f$-realizable, and
otherwise the predicate $Ar^*$-formula $\neg B(k,\ell)$ is $f$-realizable.
Consider a unary function $h(n)=\overline{\mathsf{sg}}((g(n,\tilde n)_0)$.
Obviously, $h$ is a general recursive function.
Therefore there exists a natural number $e$ such that $h(n)=\{e\}(n)$ for all~$n$.
The function $h$ has the following property:
$$h(n)=
\begin{cases}
0 &\mbox{if }\mbox{the predicate }Ar^*\mbox{-formula }\neg B(n,\tilde n)\mbox{ is }f\mbox{-realizable},\\
1 &\mbox{if }\mbox{the predicate }Ar^*\mbox{-formula }B(n,\tilde n)\mbox{ is }f\mbox{-realizable}.
\end{cases}
$$
Indeed, if the predicate $Ar^*$-formula $\neg B(n,\tilde n)$ is $f$-realizable, then $(g(n,\tilde n))_0=1$ and
$h(n)=\overline{\mathsf{sg}}((g(,n,\tilde n))_0)=0$.
If the predicate $Ar^*$-formula $B(n,\tilde n)$ is $f$-realizable, then $(g(n,\tilde n))_0=0$ and 
$h(n)=\overline{\mathsf{sg}}((g(n,\tilde n))_0)=1$.

Since the interpretation $f$ is a model of the predicate $Ar$-formula $Q$, it follows that the predicate $Ar$-formula $A_{27}$, i.e.
\begin{equation} \label{f3.70}
\forall y,z\,\neg\neg\exists v\forall x\,(x\le z\to(B(v,x)\leftrightarrow H(y,x))),
\end{equation}
is $f$-realizable.
Suppose there is a natural number $b$ that is not $f$-standard.
Since the predicate $Ar$-formula~(\ref{f3.70}) is $f$-realizable, the predicate $Ar^*$-formula
\begin{equation} \label{f3.71}
\neg\neg\exists v\forall x\,(x\le b\to(B(v,x)\leftrightarrow H(\tilde e,x)))
\end{equation}
is $f$-realizable.
Now suppose that the predicate $Ar^*$-formula
\begin{equation} \label{f141}
\exists v\forall x\,(x\le b\to(B(v,x)\leftrightarrow H(\tilde e,x)))
\end{equation}
is $f$-realizable.
Then there exists a natural number $a$ such that the predicate $Ar^*$-formula
$\forall x\,(x\le b\to(B(a,x)\leftrightarrow H(\tilde e,x)))$
is $f$-realizable.
It follows that for any $n$, the predicate $Ar^*$-formula 
$\tilde n\le b\to(B(a,\tilde n)\leftrightarrow H(\tilde e,\tilde n))$
is $f$-realizable
By Proposition~\ref{lemma8}, the predicate $Ar^*$-formula $\tilde n\le b$ is $f$-realizable.
Then for any $n$, the predicate $Ar^*$ formula $B(a,\tilde n)\leftrightarrow H(\tilde e,\tilde n))$ is $f$-realizable.
Thus for any $n$ we have:

(a) the predicate $Ar^*$-formula $B(a,\tilde n)$$f$ is $f$-realizable if and only if the predicate $Ar^*$-formula $H(\tilde e,\tilde n)$ is 
$f$-realizable.

On the other hand, by Proposition~\ref{p122}, for any $n$,

(b) the predicate $Ar^*$-formula $H(\tilde e,\tilde n)$ is $f$-realizable if and only if $h(n)=0$.

Further, by the property of the function $h$ we have:

(c) $h(n)=0$ if and only if the predicate $Ar^*$-formula $\neg B(n,\tilde n)$ is $f$-realizable.

It follows from the statements (a), (b), and (c) that the following equivalence holds for the number $a$ and any natural number $n$:
the predicate $Ar^*$-formula $B(a,\tilde n)$ is $f$-realizable if and only if the predicate $Ar^*$-formula $\neg B(n,\tilde n)$ is $f$-realizable,
and we get a contradiction if $n=a$.
This contradiction means that the predicate $Ar^*$-formula~(\ref{f141}) is not $f$-realizable.
Then its negation is $f$-realizable contrary to the $f$-realizability of the predicate $Ar^*$-formula~(\ref{f3.71}).
Thus, we have proved that every natural number is $f$-standard.~\fin

\begin{prop} \label{p121} 
Suppose an interpretation $f$ is a model of the predicate $Ar$-formula $Q$ and a number $e$ such that $e\rf Q$ is given; then for any natural number $k$, it is possible to 
effectively find a natural number $n$ such that the predicate $Ar^*$-formula $[n](k)$ is $f$-realizable,
and an $f$-realization of this formula can be found effectively.
\end{prop}

\proof
Suppose $f$ is a model of the predicate $Ar$-formula $Q$,  $e\rf Q$, and a natural number $k$ is given.
By Proposition~\ref{p118}, for any natural $n$, one can effectively find an $f$-realization of the predicate $Ar^*$-formula~(\ref{f58}) and then
an $f$-realization of the predicate $Ar^*$-formula $[n](k)\lor\neg[n](k)$, thus one can effectively check whether
the predicate $Ar^*$-formula $[n](k)$ is $f$-realizable and, if so, find its $f$-realization.
Sequentially iterating over the natural numbers starting from 0, we will find the number $n$ such that the predicate $Ar^*$-formula $[n](k)$ is 
$f$-realizable since otherwise the number $k$ would not be $f$-standard contrary to Theorem~\ref{lemma11}.~\fin

By $\nu$ denote a general recursive function which takes each $k$ to a natural number $n$ such that the predicate $Ar^*$-formula $[n](k)$ is 
$f$-realizable.
Note that $\nu(\tilde n)=n$ for any $n$ by Proposition~\ref{p61}.

\begin{prop} \label{p84}
Suppose an interpretation $f$ is a model of the predicate $Ar$-formula $Q$; then there are unary general recursive functions $\delta$ and 
$\epsilon$ such that $\delta(k)\rf E(k,\widetilde{\nu(k)})$ and $\epsilon(k)\rf E(\widetilde{\nu(k)},k)$ for any natural number~$k$.
\end{prop}

\proof
Let a natural number~$k$ be given.
Put $n=\nu(k)$; then the predicate $Ar^*$-formula $[n](k)$ is $f$-realizable and by Proposition~\ref{p121}, some its $f$-realization $p$ can be 
found effectively.
On the other hand, the predicate $Ar^*$-formula $[n](\tilde n)$ is also $f$-realizable and by Proposition~\ref{p62}, some its $f$-realization
$q$ can be found effectively from the number $n$ and, consequently, from the number~$k$.
By Proposition~\ref{p118}, one can effectively find some $f$-realization $a$ of the predicate $Ar$-formula~(\ref{f55}).
Then $\{\{\{a\}(k)\}(\nu(k))\}(2^p\cdot3^q)\rf E(k,\widetilde{\nu(k)})$, thus one can put
$\delta(k)=\{\{\{a\}(k)\}(\nu(k))\}(2^p\cdot3^q).$
Since $\forall x,y\,(x=y\to y=x)$ is deducible in $\mathsf{IRA}$, an $f$-realization of the 
predicate $Ar$-formula $\forall x,y\,(E(x,y)\to E(y,x))$ can be found effectively by Proposition~\ref{p55}, and then one can find a number $b$ 
such that $b\rf E(k,\widetilde{\nu(k)})\to E(\widetilde{\nu(k)},k)$, so we can put $\epsilon(k)=\{b\}(\delta(k)).$~\fin

\subsection{The scheme theorem} \label{ssor}
{\it The axiom of equality} for an $n$-ary predicate variable $P$ is the scheme
\begin{equation} \label{f98}
\forall \vec x\forall\vec y\,(E(x_1,y_1)\,\&\ldots\&\,E(x_n,y_n)\to(P(\vec x)\to P(\vec y))),
\end{equation}
where $\vec x$ and $\vec y$ are the lists of variables $x_1,\ldots,x_n$ and $y_1,\ldots,y_n$ respectively.

\begin{prop} \label{p98} 
For any $n$, the number $\Lambda x_1.\ldots\Lambda x_n.\Lambda y_1.\ldots\Lambda y_n.\Lambda a.\Lambda z.z$
is an $f$-realization of the scheme~(\ref{f98}) for any interpretation~$f$.
\end{prop}
\proof
This is obvious because for any $k_1,\ldots,k_n,\ell_1,\ldots,\ell_n$, if 
$$a\rf (E(k_1,\ell_1)\,\&\ldots\&\,E(k_n,l_n)),$$
i.e. 
$$a\rk(E(k_1,\ell_1)\,\&\ldots\&\,E(k_n,l_n)),$$ 
then the lists $\vec k=k_1,\ldots,k_n$ and $\vec\ell=\ell_1,\ldots,\ell_n$ coincide, and then 
the $L^*$-formulas $P(\vec k)$ and $P(\vec\ell)$ have the same $f$-realizations.~\fin

Let $Eq(E,P_1,\ldots,P_m)$ denote the conjunction of the axioms of equality for the predicate variables $P_1,\ldots,P_m$.

\begin{prop} \label{p99} 
There is a natural number $e$ such that $e$ is an $f$-realization of the scheme
$Eq(E,P_1,\ldots,P_m)$ for any interpretation~$f$.
\end{prop}

\proof
This is a direct consequence of Proposition~\ref{p98}.~\fin

Proposition~\ref{p57} can be easily extended in the following way.

\begin{prop} \label{p112}
If an interpretation $f$ is a model of the predicate formula $Q(Z,E,S,A,M)\,\&\,Eq(E,P_1,\ldots,P_m)$, then for any  $LS^*$-scheme 
$\Phi(\vec z,\vec x)$, the scheme
$\forall\vec z\forall\vec x,\vec y\,(E(\vec x,\vec y)\to(\Phi(\vec z,\vec x)\to\Phi(\vec z,\vec y)))$
is $f$-realizable.
\end{prop}

This yields the following

\begin{prop} \label{p113}
If an interpretation $f$ is a model of the predicate formula $Q(Z,E,S,A,M)\,\&\,Eq(E,P_1,\ldots,P_m)$, then for any  $LS^*$-scheme $\Phi(\vec x)$, there 
is an $(2n+2)$-place partial recursive function $\gamma_{\Phi}$ such that for any natural $d,e,\vec k,\vec\ell$, if $d\rf E(\vec k,\vec\ell)$ and 
$e\rf\Phi(\vec k)$, then $\gamma_\Phi(d,e,\vec k,\vec \ell)\rf\Phi(\vec\ell)$.
\end{prop}

\begin{prop} \label{p86} 
Suppose an interpretation $f$ is a model of the predicate formula $Q(Z,E,S,A,M)\,\&\,Eq(E,P_1,\ldots,P_m)$; then for any  $LS^*$-scheme $\Phi(x)$, 
there are binary partial recursive functions $\rho_\Phi$ and $\sigma_\Phi$ such that whatever natural numbers $a$ and $k$ are,

1) if $a\rf\Phi(k)$, then $\rho_\Phi(a,k)\rf\Phi(\widetilde{\nu(k)})$;

2) if $a\rf\Phi(\widetilde{\nu(k)})$, then $\sigma_\Phi(a,k)\rf\Phi(k)$.
\end{prop}

\proof
Let $\Phi(x)$ be an $LS^*$-scheme.
It follows from Propositions~\ref{p113} and~\ref{p84} that one can put
$$\rho_\Phi(a,k)=\gamma_\Phi(\delta(k),a,k,\widetilde{\nu(k)}),\;\sigma_\Phi(k)=\gamma_\Phi(\epsilon(k),a,k,\widetilde{\nu(k)}).$$
\fin

Let $\Phi(P_1,\ldots,P_m)$ be a closed scheme.
By $\Phi^*$ we denote the predicate formula
\begin{equation} \label{f96}
Q(Z,E,S,A,M)\,\&\,Eq(E,P_1,\ldots,P_m)\to\Phi(Z,E,S,A,M,P_1,\ldots,P_m).
\end{equation}

\begin{thm} \label{t115} 
For any closed scheme $\Phi$, if the predicate formula $\Phi^*$ is absolutely realizable, then $\Phi$ is absolutely realizable.
\end{thm}
\proof
Suppose the predicate formula~(\ref{f96}) is absolutely realizable, i.e. there exists a number $e$ such that
$$e\rf(Q(Z,E,S,A,M)\,\&\,Eq(E,P_1,\ldots,P_m)\to\Phi(P_1,\ldots,P_m))$$
for any interpretation~$f$.
In~particular, this holds for any {\it standard} interpretation $f$, where $e\rf\Psi\Leftrightarrow e\rk\Psi$ if $\Psi$ is an atomic 
$Ar^*$-formula.
Obviously, if $f$ is a standard interpretation, then for any closed predicate $Ar^*$-formula $\Psi$ and natural number $e$,
$e\rf\Psi\Leftrightarrow e\rk\Psi$.

It follows from Propositions~\ref{p15} and~\ref{p99} that in a standard interpretation $f$, the predicate formulas $Q(Z,E,S,A,M)$ and 
$Eq(E,P_1,\ldots,P_m)$ are $f$-realizable, i.e. there are numbers $q$ and $p$ such that 
$$q\rf Q(Z,E,S,A,M),\;p\rf Eq(E,P_1,\ldots,P_m)$$ 
and $q$ and $p$ do not depend on $f$, thus $\{\{e\}(2^q\cdot3^p)\rf\Phi(P_1,\ldots,P_m)$ for any interpretation~$f$.
Therefore, the scheme $\Phi$ is absolutely realizable.~\fin

The converse statement also holds, but we will have to do some preliminary work to prove it.

Suppose an interpretation $f$ is a model of the predicate $Ar$-formula $Q$ and $a\rf Q$.
By Proposition~\ref{p62}, for the number $a$ and every natural number $n$, it is possible to effectively construct the number $\tilde n$ and 
an $f$-realization of the predicate $Ar^*$-formula $[n](\tilde n)$. By Proposition~\ref{p121}, for every natural number $k$, it is 
possible to construct the number $\nu(k)$ and an $f$-realization of the predicate $Ar^*$-formula $[\nu(k)](k)$.
Consider a standard interpretation $f^\circ$, where for each (say, $m$-ary) predicate variable $P$ (other than $Z,E,S,A,M$) and for
any natural numbers $k_1,\ldots,k_m$
\begin{equation} \label{f104} 
f^\circ(P)(k_1,\ldots,k_m)=f(P)(\widetilde{k_1},\ldots,\widetilde{k_m}). 
\end{equation}

\begin{prop} \label{p1210} 
Suppose an interpretation $f$ is a model of the predicate formula $Q(Z,E,S,A,M)\,\&\,Eq(E,P_1,\ldots,P_m)$; then for any predicate formula 
$\Phi(Z,E,S,A,M,P_1,\ldots,P_n,\vec x)$, where $\vec x=x_1,\ldots,k_n$, there are $(n+1)$-place partial recursive functions $\phi_\Phi$ and $\psi_\Phi$ such that for any
$e\in\mathbb N$ and any list of natural numbers $\vec k=k_1,\ldots,k_n$, 

1) if $e\ro\Phi(P_1,\ldots,P_n,\vec k)$, then $\phi_\Phi(e,\vec k)\rf\Phi(Z,E,S,A,M,P_1,\ldots,P_n,\widetilde{\vec k})$,

2) if $e\rf\Phi(Z,E,S,A,M,P_1,\ldots,P_n,\widetilde{\vec k})$, then $\psi_\Phi(e,\vec k)\ro\Phi(P_1,\ldots,P_n,\vec k)$, where $\widetilde{\vec k}$ is the list $\widetilde{k_1},\ldots,\widetilde{k_n}$.
\end{prop}

\proof
Induction on the construction of a predicate formula~$\Phi$.

Let $\Phi(Z,E,S,A,M,P_1,\ldots,P_n,\vec x)$ be an atomic $Ar$-formula.
Denote it $\Phi(\vec x)$.
Let $\vec k$ be a list of natural numbers.
Note that $\Phi(\vec x)$ is a predicate form of a suitable arithmetical formula~$\Psi(\vec x)$.

If $e\ro\Phi(\vec k)$, then $e\rk\Phi(\vec k)$, thus the $Ar^*$-formula $\Phi(\vec k)$ is true, hence the arithmetical formula $\Psi(\vec k)$ is 
true.
Then $\mathsf{IRA}\vdash\Psi(\vec k)$ and by Proposition~\ref{p64}, one can effectively find an $f$-realization of the predicate $Ar^*$-formula
$\Phi(\widetilde{k_1},\ldots,\widetilde{k_n})$.
This $f$-realization should be taken as the value for $\phi_\Phi(e,\vec k)$.

Conversely, suppose $e\rf\Phi(\widetilde{k_1},\ldots,\widetilde{k_n})$. 
Then the $Ar^*$-formula $\Phi(\vec k)$ is true, since otherwise $\neg\Psi(\vec k)$ is true and $\mathsf{IRA}\vdash\neg\Psi(\vec k)$; 
by Proposition~\ref{p64}, the predicate $Ar^*$-formula $\neg\Phi(\widetilde{k_1},\ldots,\widetilde{k_n})$ is $f$-realizable, but this is 
impossible.
Thus the $Ar^*$-formula $\Phi(\vec k)$ is true and $0\rk\Phi(\vec k)$.
Then $0\ro\Phi(\vec k)$ and one can put $\psi_\Phi(e,\vec k)=0$.

Let $\Phi(Z,E,S,A,M,P_1,\ldots,P_n,\vec x)$ be an atomic predicate formula $P(\vec x)$, where $P$ is an $n$-ary
predicate variable.
If $e\ro P(\vec k)$, then $e\in f^\circ(P)(\vec k)$.
By~(\ref{f104}), we have
$e\in f^\circ(P)(\vec k)\Leftrightarrow e\in f(P)(\widetilde{\vec k}),$
thus one can put
$$\phi_\Phi(e,\vec k)=e;\;\psi_\Phi(e,\vec k)=e.$$

If $\Phi(Z,E,S,A,M,P_1,\ldots,P_n,\vec x)$ is of the form
$\Phi_1\,\&\,\Phi_2$
and there are functions $\phi_{\Phi_1}$, $\psi_{\Phi_1}$, $\phi_{\Phi_2}$, and $\psi_{\Phi_2}$,
then for any $e,\vec k$, one can put
$$\phi_\Phi(e,\vec k)\simeq2^{\phi_{\Phi_1}((e)_1,\vec k)}\cdot3^{\phi_{\Phi_2}((e)_2,\vec k)};\;
\psi_\Phi(e,\vec k)\simeq2^{\psi_{\Phi_1}((e)_1,\vec k)}\cdot3^{\psi_{\Phi_2}((e)_2,\vec k)}.$$

If $\Phi(Z,E,S,A,M,P_1,\ldots,P_n,\vec x)$ is of the form
$\Phi_1\lor\Phi_2,$
then for any $e,\vec k$, one can put
$$\phi_\Phi(e,\vec k)\simeq
\begin{cases}
2^0\cdot3^{\phi_{\Phi_1}((e)_1,\vec k)}&\mbox{if }(e)_0=0,\\
2^1\cdot3^{\phi_{\Phi_2}((e)_1,\vec k)}&\mbox{if }(e)_0\not=0;
\end{cases}
$$
$$\psi_\Phi(e,\vec k)\simeq
\begin{cases}
2^0\cdot3^{\psi_{\Phi_1}((e)_1,\vec k)}&\mbox{if }(e)_0=0,\\
2^1\cdot3^{\psi_{\Phi_2}((e)_1,\vec k)}&\mbox{if }(e)_0\not=0.
\end{cases}
$$

If $\Phi(Z,E,S,A,M,P_1,\ldots,P_n,\vec x)$ is of the form
$\Phi_1\to\Phi_2,$
then for any $e,\vec k$, one can put
$$\phi_\Phi(e,\vec k)=\Lambda x.\phi_{\Phi_2}(\{e\}(\psi_{\Phi_1}(x,\vec k)),\vec k);\;
\psi_\Phi(e,\vec k)=\Lambda x.\psi_{\Phi_2}(\{e\}(\phi_{\Phi_1}(x,\vec k)),\vec k).$$

If $\Phi(Z,E,S,A,M,P_1,\ldots,P_n,\vec x)$ is of the form
$\neg\Psi,$
then for any $e,\vec k$, one can put
$\phi_\Phi(e,\vec k)=\Lambda x.0,\;\psi_\Phi(e,\vec k)=\Lambda x.0.$

Suppose $\Phi(Z,E,S,A,M,P_1,\ldots,P_n,\vec x)$ is of the form
$$\forall x\,\Psi(Z,E,S,A,M,P_1,\ldots,P_n,x,\vec x)$$
and there are $(n+2)$-place partial recursive functions $\phi_\Psi$ and $\psi_\Psi$ such that for any $e,k,\vec k$, if
$e\ro\Psi(k,\vec k)$, then $\phi_\Psi(e,k,\vec k)\rf\Psi(\tilde k,\widetilde{\vec k})$, and
if $e\rf\Psi(\tilde k,\widetilde{\vec k})$, then $\psi_\Psi(e,k,\vec k)\ro\Psi(k,\vec k)$.
In this case, for any $e,\vec k$, one can put
$$\phi_\Phi(e,\vec k)=\Lambda x.\sigma_\Theta(\phi_\Psi(\{e\}(\nu(x)),\nu(x),\vec k),x),$$
where $\Theta(x)$ is the predicate $Ar^*$-formula $\Psi(x,\widetilde{k_1},\ldots,\widetilde{k_n})$ and $\sigma_\Theta$ is the function from
Proposition~\ref{p86};
$$\psi_\Phi(e,\vec k)=\Lambda x.\psi_\Psi(\{e\}(\tilde x),x,\vec k).$$

If the predicate formula $\Phi(Z,E,S,A,M,P_1,\ldots,P_n,\vec x)$ is of the form
$$\exists x\,\Psi(Z,E,S,A,M,P_1,\ldots,P_n,\vec x),$$
then for any $e,\vec k$ one can put
$$\phi_\Phi(e,\vec k)\simeq2^{\widetilde{(e)_0}}\cdot3^{\phi_\Psi((e)_1,(e)_0,\vec k)};$$
$$\psi_\Phi(e,\vec k)\simeq2^{\nu((e)_0)}\cdot3^{\psi_\Psi(\rho_\Theta((e)_1,(e)_0),\nu((e)_0),\vec k)},$$
where $\rho_\Theta$ is the function from Proposition~\ref{p86}.~\fin

\begin{thm} \label{t129} 
If the scheme $\Phi$ is absolutely realizable, then the predicate formula $\Phi^*$ is absolutely realizable.
\end{thm} 
                                                        
\proof
Suppose the scheme $\Phi$ is absolutely realizable.
We prove that the predicate formula~(\ref{f96}) is absolutely realizable.
To do this, we describe an algorithm that for any interpretation $f$, constructs an $f$-realization of the predicate formula 
$\Phi(Z,E,S,A,M,P_1,\ldots,P_m)$ if an $f$-realization of the predicate formula~$Q\,\&\,Eq(E,P_1,\ldots,P_m)$ is given.
Suppose an $f$-realization of the predicate formula~$Q\,\&\,Eq(E,P_1,\ldots,P_m)$ is given.
It is of the form $2^q\cdot3^p$, where $q\rf Q$, $p\rf Eq(E,P_1,\ldots,P_m)$.
This means, in particular, that the interpretation $f$ is a model of the predicate $Ar$-formula $Q$.
By Proposition~\ref{p62}, for the number $q$ and every natural number $n$ one can effectively construct the number $\tilde n$.
Consider the standard interpretation $f^\circ$ defined by~(\ref{f104}).
Then $a\ro\Phi$ and by Proposition~\ref{p1210}, we have $\phi_\Phi(a)\rf\Phi$.
Thus for any $f$realization of $Q\,\&\,Eq(E,P_1,\ldots,P_m)$ we can effectively find an $f$-realization of the formula $\Phi$, and
the described algorithm does not depend on the interpretation $f$ itself.
Therefore, the predicate formula $\Phi^*$ is absolutely realizable.~\fin

\begin{thm} \label{main} 
For any closed scheme $\Phi$
, the predicate formula $\Phi^*$ is absolutely realizable if and only if $\Phi$ is absolutely
realizable.
\end{thm}

\proof
This is a direct consequence of Theorems~\ref{t115} and~\ref{t129}.~\fin

\section{Constructive predicate calculus} \label{skip} 

\subsection{Extrended Church thesis} \label{ssotch}

Recall that $G(x,y,z)$ is a $\Sigma$-formula expressing the predicate $\{x\}(y)=z$.
By Proposition~\ref{p69}, there exists a triple partial recursive function $\psi$ (namely, $\alpha_G$) such that whatever
natural numbers $k,\ell,m$ are, if $\{k\}(\ell)=m$, then $!\psi(k,\ell,m)$ and $\psi(k,\ell,m)\rk G(k,\ell,m)$.
The scheme $nCT$ was introduced in Subsection~\ref{sspf}.

\begin{thm} \label{t21}
The scheme $nCT$ is absolutely realizable.
\end{thm}

\proof
Let $e=\Lambda a.2^{l(a)}\cdot3^{r(a)}$, where
$$l(a)=\Lambda x.(g(a,x))_0,\;g(a,x)=\{\{a\}(x)\}(0),$$ 
$$r(a)=\Lambda x.\Lambda u.2^{(g(a,x))_0}\cdot3^{2^d\cdot3^{(g(a,x))_1}},\;d=\psi(l(a),x,(g(a,x))_0).$$
We prove that $e$ $f$-realizes the scheme $nCT$ for any interpretation~$f$.
This means that for any natural $a$, if
\begin{equation} \label{f22}
a\rf\forall x\,(\neg P(x)\to\exists y\,Q(x,y)), 
\end{equation}
then $2^{l(a)}\cdot3^{r(a)}\rf\exists z\forall x\,(\neg P(x)\to\exists y\,(G(z,x,y)\,\&\,Q(x,y))),$
i.e. for any $k$, 
$\Lambda u.2^{g(a,k))_0}\cdot3^{2^d\cdot3^{g(a,k))_1}}\rf(\neg P(k)\to\exists y\,(G(l(a),k,y)\,\&\,Q(k,y))),$
in other words, whatever a number $b$ is, if $b\rf\neg P(k)$, then
$$2^d\cdot3^{(g(a,k))_1}\rf G(l(a),k,(g(a,k))_0)\,\&\,Q(k,(g(a,k))_0)).$$
This means that
\begin{equation} \label{f23}
d\rf G(l(a),k,(g(a,k))_0),
\end{equation}
\begin{equation} \label{f235}
(g(a,k))_1\rf Q(k,(g(a,k))_0).
\end{equation}
So, assume that $b\rf\neg P(k)$ for some~$b$.
Then by Proposition~\ref{p5}, $0\rf\neg P(k)$.
This and~(\ref{f22}) imply that the value $g(a,k)$ is defined and $g(a,k)\rf\exists y\,Q(k,y)$, i.e.~(\ref{f235}) holds.
The condition~(\ref{f23}) also holds because $\{l(a)\}(k)=(g(a,k))_0$.~\fin

\subsection{Calculus $\mathsf{MQC}$} \label{ssimqc}

By means of $\mathsf{MQC}$ we denote the calculus obtained by adding to the axioms of $\mathsf{IQC}$ the Markov principle $M$ and the formula~$nCT^*$.
This means that all substitutional instances of the predicate formulas $M$ and $nCT^*$ are axioms of~$\mathsf{MQC}$.

\begin{thm} \label{t132} 
Only absolutely realizable predicate formulas are deduced in the calculus $\mathsf{MQC}$.
\end{thm}

\proof
This is a direct consequence of Theorems~\ref{t32}, \ref{t18}, \ref{t21},~\ref{t129}.~\fin

By means of $CT(Q)$ we denote the scheme
$$\forall x\,\exists y\,Q(x,y)\to\exists z\,\forall x\,\exists y\,(G(z,x,y)\,\&\,Q(x,y)).$$

\begin{prop} \label{p133} 
The predicate formula $CT^*(Q)$ is deducible in the calculus $\mathsf{MQC}$.
\end{prop}

\proof
Consider the following substitutional instance of the predicate formula $nCT^*(P,Q)$: substitute for $P$ a formula $\Phi(x)$ such that $\neg\Phi(x)$ is deducible in $\mathsf{IQC}$.
For example, one can take $Z(x)\,\&\,\neg Z(x)$ as $\Phi(x)$.
Then it is obvious that $CT^*(Q)$ is derived from $nCT^*(\Phi(x),Q)$ in $\mathsf{IQC}$, hence in $\mathsf{MQC}$.~\fin

By means of $\mathsf{CA}$ we denote an arithmetical theory based on the system of Peano axioms and the calculus
$\mathsf{MQC}$.
It is obvious that every formula deducible in intuitionistic arithmetic $\mathsf{HA}$ is also deducible in~$\mathsf{CA}$.

\begin{thm} \label{t133} 
Every formula derived in Markov arithmetic $\mathsf{MA}$ is deducible in~$\mathsf{CA}$.
\end{thm}

\proof
It is enough to prove that any arithmetical formula $\Theta$ of the form
$$\forall x\,(\neg\Phi(x)\to\exists y\,\Psi(x,y))\to\exists z\forall x\,(\neg\Phi(x)\to\exists y\,(G(z,x,y)\,\&\,\Psi(x,y)))$$
obtained by the scheme nCT is deducible in~$\mathsf{CA}$.
Note that the arithmetical formula $\Theta^*$, i.e. $Q\,\&\,Eq(\Phi,\Psi)\to\Theta$, is an axiom of the theory~$\mathsf{CA}$.
Now recall that $Q$ is the conjunction of the axioms of Robinson arithmetic derived in $\mathsf{HA}$, hence in
$\mathsf{CA}$, and the formulas $A_{26}$, $A_{27}$, and $A_{28}$, whose deducibility in $\mathsf{HA}$ was established above
(see respectively~(\ref{f33}), Proposition~\ref{p13}, (\ref{f21})).
Hence $Q$ is deducible in~$\mathsf{CA}$.
The formula $Eq(\Phi,\Psi)$ is also deducible.
It follows that $\Theta$ is deducible in $\mathsf{CA}$, as was to be proved.~\fin

It is of interest if the converse is true, that every arithmetical formula deducible in $\mathsf{CA}$ is deducible in Markov 
arithmetic~$\mathsf{MA}$.
Obviously, to do this, we need to prove that every arithmetical substitutional instance of the predicate formula $nCT^*$ is derived in Markov 
arithmetic $\mathsf{MA}$.
The study of this question was not the purpose of this article.

\begin{thm} \label{t134} 
There is a predicate formula deducible in the calculus $\mathsf{MAC}$, but not deducible in the classical predicate calculus~$\mathsf{CQC}$.
\end{thm}

\proof
Let $\Phi(x)$ be the arithmetical formula $\exists y\,(Z(y)\,\&\,G(x,x,y))$.
Consider the predicate $Ar$-formula $\neg\forall x\,(\neg\Phi(x)\lor\Phi(x))$.
Denote it~$\Psi$.
It is quite obvious that the predicate formula $\Psi^*$ is not classically valid, therefore, is not deducible in the classical
predicate calculus $\mathsf{CQC}$.
We prove that the predicate formula $\Psi^*$ is deducible in $\mathsf{MQC}$.
The formula $\Psi^*$ is of the form $Q\to\Psi$, so it is enough to deduce in $\mathsf{MQC}$ the formula $\Psi$ from the 
hypothesis~$Q$.
Let's prove that in $\mathsf{IQC}$, the formula
$$(\neg\Phi(x)\lor\Phi(x))\leftrightarrow\exists y\,((Z(y)\to\neg\Phi(x))\,\&\,(\neg Z(y)\to\Phi(x))$$
is derived from $Q$.
First we prove that
$$Q,\,\neg\Phi(x)\lor\Phi(x)\vdash\exists y\,((Z(y)\to\neg\Phi(x))\,\&\,(\neg Z(y)\to\Phi(x)).$$
It is sufficient to prove that
\begin{equation} \label{f137} 
Q,\,\neg\Phi(x)\vdash\exists y\,((Z(y)\to\neg\Phi(x))\,\&\,(\neg Z(y)\to\Phi(x)) 
\end{equation}
and
\begin{equation} \label{f138} 
Q,\,\Phi(x)\vdash\exists y\,((Z(y)\to\neg\Phi(x))\,\&\,(\neg Z(y)\to\Phi(x)).
\end{equation}
Obviously, $Q,\,\neg\Phi(x),\,Z(y)\vdash((Z(y)\to\neg\Phi(x))\,\&\,(\neg Z(y)\to\Phi(x)))$, thus
$$Q,\,\neg\Phi(x),\,\exists y\,Z(y)\vdash\exists y\,((Z(y)\to\neg\Phi(x))\,\&\,(\neg Z(y)\to\Phi(x))).$$
It remains to note that $\exists y\,Z(y)$ is the formula $A_{14}$, which is a conjunctive member of $Q$ and is therefore deducible 
from it.
Thus~(\ref{f137}) is proved.
To prove~(\ref{f138}) note that
$$Q,\,\Phi(x),\,\neg Z(y)\vdash((Z(y)\to\neg\Phi(x))\,\&\,(\neg Z(y)\to\Phi(x)))$$
and $\exists y\,\neg Z(y)$ is derived from $A_{14}$, $A_{17}$, and~$A_2$.

The formula $(Z(y)\to\neg\Phi(x))\,\&\,(\neg Z(y)\to\Phi(x)$ is denoted by~$\Theta(x,y)$.
Thus it is sufficient to deduce $\neg\forall x\exists y\,\Theta(x,y)$ from the hypothesis ~$Q$ in $\mathsf{IQC}$.

By Proposition~\ref{p133}, the predicate formula $CT^*(\Theta)$, i.e.
$$Q\,\&\,Eq(E(x,y),\Theta(x,y))\to(\forall x\,\exists y\,\Theta(x,y)\to\exists z\,\forall x\,\exists y\,(G(z,x,y)\,\&\,\Theta(x,y))),$$
is deduced in $\mathsf{MQC}$.
The premise of this formula is deducible from the hypothesis $Q$, thus the conclusion 
$$\forall x\,\exists y\,\Theta(x,y)\to\exists z\,\forall x\,\exists y\,(G(z,x,y)\,\&\,\Theta(x,y))$$
is also deducible.
To prove the deducibility of $\Psi$ in $\mathsf{MQC}$, it is enough to deduce a contradiction from the hypothesis 
$\exists y\,\Theta(x,y)$.
We see that the formula 
$$\exists z\,\forall x\,\exists y\,(G(z,x,y)\,\&\,\Theta(x,y)))$$
is derived from this hypothesis.
It is enough to deduce the contradiction from $\forall x\,\exists y\,(G(z,x,y)\,\&\,\Theta(x,y))$.
Note that $\exists y\,(G(z,z,y)\,\&\,\Theta(z,y))$, i.e.
$$\exists y\,(G(z,z,y)\,\&\,((Z(y)\to\neg\exists v\,(Z(v)\,\&\,G(z,z,v))\,\&$$
$$\&\,(\neg Z(y)\to\exists v\,(Z(v)\,\&\,G(z,z,v))),$$
is derived from this hypothesis
It is enough to deduce a contradiction from the hypothesis
$$G(z,z,y)\,\&\,((Z(y)\to\neg\exists v\,(Z(v)\,\&\,G(z,z,v))\,\&$$
$$\&\,(\neg Z(y)\to\exists v\,(Z(v)\,\&\,G(z,z,v)).$$
The formulas $G(z,z,y)$,
$$Z(y)\to\neg\exists v\,(Z(v)\,\&\,G(z,z,v)),$$
$$\neg Z(y)\to\exists v\,(Z(v)\,\&\,G(z,z,v))$$
are derived from this hypothesis.
The formula $Z(y)\lor\neg Z(y)$ is derived from $Q$, so it is enough to deduce a contradiction from each of the hypotheses $Z(y)$ and 
$\neg Z(y)$.
The formulas $\neg\exists v\,(Z(v)\,\&\,G(z,z,v))$ and $Z(y)\,\&\,G(z,z,y)$ are derived from the hypothesis $Z(y)$, thus a contradiction follows.
The formula $\exists v\,(Z(v)\,\&\,G(z,z,v))$ is derived fFrom the hypothesis $\neg Z(y)$.
It remains to noted that $E(y,v)$ is derived from the hypotheses $Z(v)\,\&\,G(z,z,v)$ and $G(z,z,y)$ and $A_{26}$, and 
further we get $Z(y)$ from $A_{16}$, and this leads to a contradiction.
Thus the deducibility of $\Psi^*$ is proved.~\fin

\end{document}